# Boundary behaviour of the solution of the heat equation on the half line via the Fokas unified transform method

Andreas Chatziafratis

National and Kapodistrian University of Athens, Greece

e-mail: chatziafrati@math.uoa.gr

**Abstract.** We consider the Fokas method expression for the solution of the heat equation $u_t = u_{xx}$ on the half line with Dirichlet data and we study in detail its boundary behaviour, i.e., as $x \to 0^+$ or $t \to 0^+$ (including the case $(x,t) \to (0,0)$) by analyzing the integrals involved. We also study the boundary behaviour of the derivatives of the solution. In particular we give conditions on the data which guarantee the extension of the solution to a $C^\infty$ function on $\{x \geq 0, t \geq 0\}$.

## 1. Introduction

Our starting point is Fokas method of solving initial and boundary value problems and more specifically the form of the solution this method gives for the heat equation on the half line. (See [2,3,4,5,6,7,8,9,10,11,12,15] ). We use this equation as an illustrative example and we make a detailed study of the boundary behaviour of the solution given by the Fokas method.
Considering the equation
$$u_t = u_{xx}, \text{ for } x > 0 \text{ and } t > 0,$$
with the initial and boundary conditions
$$u(x,0) = u_0(x) \text{ for } x \geq 0 \text{ and } u(0,t) = g_0(t) \text{ for } t \geq 0,$$
the Fokas method gives the following integral representation: For $x > 0$ and $t > 0$,
$$u(x,t) = \frac{1}{2\pi} \int_{\lambda=-\infty}^{\infty} e^{i\lambda x - \lambda^2 t} \hat{u}_0(\lambda) d\lambda - \frac{1}{2\pi} \int_{\lambda \in \Gamma} e^{i\lambda x - \lambda^2 t} \hat{u}_0(-\lambda) d\lambda - \frac{i}{\pi} \int_{\lambda \in \Gamma} e^{i\lambda x - \lambda^2 t} \lambda \tilde{g}_0(\lambda, t) d\lambda \qquad (1.1)$$
where
$$\hat{u}_0(\lambda) = \int_{x=0}^{\infty} u_0(x) e^{-i\lambda x} dx, \text{ defined for } \lambda \in \mathbb{C} \text{ with } \text{Im } \lambda \leq 0,$$
$$\tilde{g}_0(\lambda, t) = \int_{\tau=0}^{t} e^{\lambda^2 \tau} g_0(\tau) d\tau, \text{ defined for } \lambda \in \mathbb{C},$$
and the contour $\Gamma$ is the oriented boundary of the domain $\Omega^- := \{\lambda \in \mathbb{C} : \text{Im } \lambda \geq 0 \text{ and } \text{Re}(\lambda^2) \leq 0\}$.

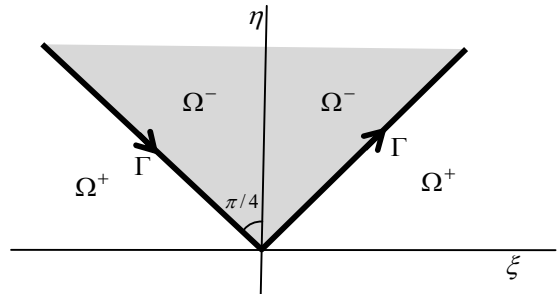

*The contour* $\Gamma = \partial \Omega^- = \{\lambda \in \mathbb{C} : \text{Im } \lambda \geq 0 \text{ and } \text{Re}(\lambda^2) = 0\}$

Date: January 2020

**Acknowledgements**
The work contained in this preprint (January 2020) is part of my M.Sc. thesis (2018-19) at the National and Kapodistrian University of Athens, Greece. I wish to express my gratitude to my Professors: N. Alikakos, G. Barbatis, T. Hatziafratis, I. Stratis, for insightful discussions and continuous support. The Onassis Foundation is also gratefully acknowledged for partial funding during that stage. The results presented herein have been published as: A. Chatziafratis, D. Mantzavinos, Boundary behavior for the heat equation on the half-line, *Math. Methods Appl. Sci.* 45, 7364-93 (2022).

Furthermore
$$\frac{\partial^n u(x,t)}{\partial x^n} = \frac{i^n}{2\pi} \int_{\lambda=-\infty}^{\infty} \lambda^n e^{i\lambda x - \lambda^2 t} \hat{u}_0(\lambda) d\lambda - \frac{i^n}{2\pi} \int_{\lambda \in \Gamma} \lambda^n e^{i\lambda x - \lambda^2 t} \hat{u}_0(-\lambda) d\lambda$$
$$- \frac{i^{n+1}}{\pi} \int_{\lambda \in \Gamma} \lambda^n e^{i\lambda x - \lambda^2 t} \lambda \tilde{g}_0(\lambda, t) d\lambda, \quad (1.2)$$

for $x > 0$, $t > 0$ and $n \in \mathbb{N}$.

Also
$$\frac{\partial^n u(x,t)}{\partial t^n} = \frac{(-1)^n}{2\pi} \int_{\lambda=-\infty}^{\infty} \lambda^{2n} e^{i\lambda x - \lambda^2 t} \hat{u}_0(\lambda) d\lambda - \frac{(-1)^n}{2\pi} \int_{\lambda \in \Gamma} \lambda^{2n} e^{i\lambda x - \lambda^2 t} \hat{u}_0(-\lambda) d\lambda$$
$$- (-1)^n \frac{i}{\pi} \int_{\lambda \in \Gamma} \lambda^{2n} e^{i\lambda x - \lambda^2 t} \lambda \tilde{g}_0(\lambda, T) d\lambda, \quad (1.3)$$

for $x > 0$, $0 < t < T$ and $n \in \mathbb{N} \cup \{0\}$.

In this paper we will study the limits of the function $u(x,t)$ (originally defined for $x > 0$ and $t > 0$ by (1.1)) and its derivatives $\partial^{k+\ell} u(x,t) / \partial x^k \partial t^\ell$, as $(x,t)$ approaches a point of the boundary of $Q := \{(x,t) \in \mathbb{R}^2 : x > 0 \text{ and } t > 0\}$. More precisely we will prove the following theorems.

**Theorem 1.1** *Given $u_0(x) \in \mathcal{S}([0,\infty))$ and $g_0(t) \in C^\infty([0,\infty))$, the function $u(x,t)$ defined by (1.1) is $C^\infty$ for $(x,t) \in Q$ and satisfies the following:*

*1st The differential equation $u_t = u_{xx}$ for $x > 0$ and $t > 0$.*

*2nd The limit conditions $\lim_{t \to 0^+} u(x,t) = u_0(x)$ (for each fixed $x > 0$) and $\lim_{x \to 0^+} u(x,t) = g_0(t)$ (for each fixed $t > 0$).*

*3rd $u(x,t) \in C^\infty([0,\infty))$ with respect to $x$ (for each fixed $t > 0$), and the functions $g_n(t)$, $t \in (0,\infty)$, defined by*

$$g_n(t) = \lim_{x \to 0^+} \frac{\partial^n u(x,t)}{\partial x^n}, \quad n \in \mathbb{N},$$

*are $C^\infty$ for $t \in (0,\infty)$.*

*4th $u(x,t) \in C^\infty([0,\infty))$ with respect to $t$ (for each fixed $x > 0$), and the functions $u_n(x)$, $x \in (0,\infty)$, defined by*

$$u_n(x) = \lim_{t \to 0^+} \frac{\partial^n u(x,t)}{\partial t^n}, \quad n \in \mathbb{N},$$

*are $C^\infty$ for $x \in (0,\infty)$.*

*5th The function $u(x,t)$ is rapidly decreasing as $x \to +\infty$, uniformly for $t$ in compact subsets of $(0,\infty)$.*

**Theorem 1.2** *If $u_0(x) \in \mathcal{S}([0,\infty))$ and $g_0(t) \in C^\infty([0,\infty))$ then the function $u(x,t)$ defined by (1.1) for $(x,t) \in Q$ extends to a $C^\infty$ function on $\overline{Q} - \{(0,0)\}$, i.e., all the derivatives*

$$\frac{\partial^{k+\ell} u(x,t)}{\partial x^k \partial t^\ell}, \quad (x,t) \in Q, \; k, \ell \in \mathbb{N} \cup \{0\},$$

*extend continuously to $\{(x,t) \in \mathbb{R}^2 : x \geq 0 \text{ and } t \geq 0\} - \{(0,0)\}$.*

In order to have limits of the function $u(x,t)$ (or certain of its derivatives) also when $(x,t) \to (0,0)$, we should make some further assumptions on the functions $u_0(x)$ and $g_0(t)$. For example we will prove the following theorem.

**Theorem 1.3** *Suppose* $u_0(x) \in \mathcal{S}([0,\infty))$ *and* $g_0(t) \in C^\infty([0,\infty))$. *Then the functions* $u(x,t)$, $\partial u(x,t)/\partial x$, $\partial^2 u(x,t)/\partial x^2$ *and* $\partial^3 u(x,t)/\partial x^3$ *(originally defined by (1.1) and (1.2) for* $(x,t) \in Q$ *and extended to* $\overline{Q} - \{(0,0)\}$ *by Theorem 1.2) satisfy the following:*

*1st If, in addition,* $u_0(0) = g_0(0)$ *then*

$$\lim_{\overline{Q} \ni (x,t) \to (0,0)} u(x,t) = u_0(0) \quad \text{and} \quad \lim_{\overline{Q} \ni (x,t) \to (0,0)} \frac{\partial u(x,t)}{\partial x} = \frac{du_0(x)}{dx}\bigg|_{x=0}.$$

*2nd If* $u_0(0) = g_0(0)$ *and* $\dfrac{d^2 u_0(x)}{dx^2}\bigg|_{x=0} = \dfrac{dg_0(t)}{dt}\bigg|_{t=0}$ *then*

$$\lim_{\overline{Q} \ni (x,t) \to (0,0)} \frac{\partial^2 u(x,t)}{\partial x^2} = \frac{d^2 u_0(x)}{dx^2}\bigg|_{x=0} \quad \text{and} \quad \lim_{\overline{Q} \ni (x,t) \to (0,0)} \frac{\partial^3 u(x,t)}{\partial x^3} = \frac{d^3 u_0(x)}{dx^3}\bigg|_{x=0}.$$

Part of the difficulty in dealing with these questions is due to the fact that some of the integrals which are involved in the representation of the solution or of its derivatives do not converge absolutely, when $x = 0$ or $t = 0$, and they have to be interpreted in a generalized sense. For example the first integral in the RHS of (1.2), for $t = 0$, becomes

$$\left[ \int_{\lambda = -\infty}^{\infty} \lambda^n e^{i\lambda x - \lambda^2 t} \hat{u}_0(\lambda) d\lambda \right]_{t=0} = \int_{\lambda = -\infty}^{\infty} \lambda^n e^{i\lambda x} \hat{u}_0(\lambda) d\lambda$$

and $\left|\lambda^n e^{i\lambda x} \hat{u}_0(\lambda)\right| \cong |\lambda|^{n-1}$ as $\lambda \to \pm\infty$ ($\lambda \in \mathbb{R}$).

The second integral in the RHS of (1.3) becomes

$$\left[ \int_{\lambda \in \Gamma} \lambda^{2n} e^{i\lambda x - \lambda^2 t} \hat{u}_0(-\lambda) d\lambda \right]_{x=0} = \int_{\lambda \in \Gamma} \lambda^{2n} e^{-\lambda^2 t} \hat{u}_0(-\lambda) d\lambda$$

and $\left|\lambda^{2n} e^{-\lambda^2 t} \hat{u}_0(-\lambda)\right| \cong |\lambda|^{2n-1}$ as $\lambda \to \infty$ with $\lambda \in \Gamma$, since $\left|e^{-\lambda^2 t}\right| = 1$ for $\lambda \in \Gamma$. (See also (2.5) below.)

Similarly the third integral in the RHS of (1.2) becomes

$$\left[ \int_{\lambda \in \Gamma} \lambda^n e^{i\lambda x - \lambda^2 t} \lambda \tilde{g}_0(\lambda, t) d\lambda \right]_{x=0} = \int_{\lambda \in \Gamma} \lambda^n e^{-\lambda^2 t} \lambda \tilde{g}_0(\lambda, t) d\lambda$$

and $\left|\lambda^n e^{-\lambda^2 t} \lambda \tilde{g}_0(\lambda, t)\right| \cong |\lambda|^{n-1}$ as $\lambda \to \infty$ with $\lambda \in \Gamma$. (This follows from (2.9) below.)

The situation becomes worse when both $x = 0$ and $t = 0$, i.e., when $(x,t) \to (0,0)$. In this case more assumptions have to be imposed (on the data $u_0$ and $g_0$) in order to study the limits of the various terms. (See Theorem 1.3 and, its generalization, Theorem 8.1.)

We will be using the following – rather standard – definition of the spaces $C^\infty([0,\infty))$ and $\mathcal{S}([0,\infty))$.

*Definitions* As usual, we will say that a $C^\infty$ function $f:(0,\infty)\to\mathbb{C}$ belongs to the space $C^\infty([0,\infty))$ and we will write $f\in C^\infty([0,\infty))$ if the limit $\lim_{x\to 0^+}\dfrac{d^N f(x)}{dx^N}$ exists for every nonnegative integer $N$. Equivalently, $f\in C^\infty([0,\infty))$ if it extends to a $C^\infty$ function $\widetilde{f}:\mathbb{R}\to\mathbb{C}$.

We will say that a function $f\in\mathcal{S}([0,\infty))$ if $f\in C^\infty([0,\infty))$ and it is rapidly decreasing as $x\to+\infty$ in the sence that

$$\sup_{x\geq 0}\left|x^M \frac{d^N f(x)}{dx^N}\right|<\infty \text{ for every nonnegative integers } M \text{ and } N.$$

Equivalently, $f\in\mathcal{S}([0,\infty))$ if $f\in C^\infty([0,\infty))$ and

$$\lim_{x\to+\infty}\left|x^M \frac{d^N f(x)}{dx^N}\right|=0 \text{ for every nonnegative integers } M \text{ and } N.$$

A $C^\infty$ function $F:Q\to\mathbb{C}$, $F=F(x,t)$, $(x,t)\in Q$, is said to be rapidly decreasing as $x\to+\infty$, uniformly for $t$ in compact subsets of $(0,\infty)$, if for every $\beta>\alpha>0$,

$$\sup_{x\geq 0,\alpha\leq t\leq \beta}\left|x^M \frac{\partial^N F(x,t)}{\partial x^N}\right|<\infty \text{ for every nonnegative integers } M \text{ and } N.$$

Equivalently, $F=F(x,t)$ is rapidly decreasing as $x\to+\infty$, uniformly for $t$ in compact subsets of $(0,\infty)$, if, for every nonnegative integers $M$ and $N$,

$$\lim_{x\to+\infty}\left|x^M \frac{\partial^N F(x,t)}{\partial x^N}\right|=0, \text{ uniformly for } t \text{ in compact subsets of } (0,\infty).$$

The proof of Theorem 1.1 will be completed in section 5 and the proof of Theorem 1.2 will be given in section 7. In sections 6 and 7 we study in detail the boundary behaviour of the solution and its derivatives needed for the proof of Theorem 1.2. Finally in section 8 we give the proof of Theorem 1.3.

## 2. Preliminaries

### 2.1. Fourier's inversion formula.

If a function $f:\mathbb{R}\to\mathbb{C}$ is $L^1$, i.e., Lebesgue measurable and $\int_{-\infty}^{\infty}|f(x)|dx<\infty$, then its Fourier transform is defined by the formula

$$\hat{f}(\lambda)=\int_{x=-\infty}^{\infty}e^{-i\lambda x}f(x)dx, \text{ for } \lambda\in\mathbb{R},$$

and the resulting function $\hat{f}:\mathbb{R}\to\mathbb{C}$ is continuous and bounded, and, moreover, $\hat{f}(\lambda)\to 0$, as $\lambda\to\pm\infty$. According to the inversion formula, if the function $\hat{f}:\mathbb{R}\to\mathbb{C}$ is also $L^1$, then

$$f(x)=\frac{1}{2\pi}\int_{\lambda=-\infty}^{\infty}e^{i\lambda x}\hat{f}(\lambda)d\lambda, \text{ a.e. for } x\in\mathbb{R}.$$

For convenience we also state two basic inversion formulas which we will use in this paper and which can be applied in cases with $\hat{f}\notin L^1(\mathbb{R})$, provided that some other conditions are satisfied. More precisely, even if $\hat{f}\notin L^1(\mathbb{R})$, we can recover the original function from its Fourier transform in the following cases:

1. If $f \in L^1(\mathbb{R})$ and is piecewise continuous then

$$\frac{1}{2\pi} \lim_{\varepsilon \to 0^+} \int_{\lambda=-\infty}^{\infty} e^{-\varepsilon\lambda^2} e^{i\lambda x} \hat{f}(\lambda) d\lambda = \frac{1}{2}[f(x^+) + f(x^-)] \text{ for every } x \in \mathbb{R}, \qquad (2.1)$$

where $f(x^\pm) = \lim_{y \to x^\pm} f(y)$.

2. If the function $f : \mathbb{R} \to \mathbb{C}$ is $L^1$ and is piecewise $C^1$, then

$$\frac{1}{2\pi} \lim_{A \to \infty} \int_{\lambda=-A}^{A} e^{i\lambda x} \hat{f}(\lambda) d\lambda = \frac{1}{2}[f(x^+) + f(x^-)] \text{ for every } x \in \mathbb{R}. \qquad (2.2)$$

(See also [13, section 7.2].)

## *2.2. Lebesgue's dominated convergence theorem.*

In several instances, we will have to deal with limit processes and, in particular, with problems of whether a specific limit process can be interchanged with an integration process. Our basic tool in dealing with such problems will be Lebesgue's dominated convergence theorem. The following two versions of this theorem are the ones which we will need.

1. A condition which guarantees the validity of equation

$$\lim_{\lambda \to \lambda_0} \int_{x \in X} f(x,\lambda) d\mu(x) = \int_{x \in X} \lim_{\lambda \to \lambda_0} f(x,\lambda_0) d\mu(x)$$

is the following: It suffices that there exist a function $g(x)$ such that $|f(x,\lambda)| \leq g(x)$ for $x \in X$ and $\lambda$ close to $\lambda_0$, and with $\int_{x \in X} g(x) d\mu(x) < \infty$. (Here $d\mu$ is supposed to be an appropriate measure in the space $X$, $\lambda$ is a parameter, and the functions $f(x,\lambda)$ and $g(x)$ are supposed to be appropriately measurable functions for $x \in X$.)

2. A condition which guarantees the validity of the equation

$$\frac{\partial}{\partial \lambda}\bigg|_{\lambda=\lambda_0} \left( \int_{x \in X} f(x,\lambda) d\mu(x) \right) = \int_{x \in X} \left[ \frac{\partial f(x,\lambda)}{\partial \lambda} \bigg|_{\lambda=\lambda_0} \right] d\mu(x) \qquad (2.3)$$

is the following: It suffices that there exist a function $g(x)$ such that $|\partial f(x,\lambda)/\partial\lambda| \leq g(x)$, for $x \in X$ and for $\lambda$ close to $\lambda_0$, and with $\int_{x \in X} g(x) dx < \infty$. (Here $\lambda$ is supposed to be a real parameter and the function $f(x,\lambda)$ is assumed to be $C^1$ with respect to $\lambda$, for $\lambda$ close to $\lambda_0$ and for every $x \in X$.)

## *2.3. Jordan's type lemmas.*

The proof of these lemmas is based on the following elementary inequality:

$$\text{For } \kappa > 0, \ \int_{\theta=0}^{\pi} e^{-\kappa \sin\theta} d\theta \leq \frac{\pi}{\kappa}.$$

Also the above integral is equal to each of the following ones:

$$2\int_{\theta=0}^{\pi/2} e^{-\kappa\sin\theta} d\theta = \int_{\theta=-\pi/2}^{\pi/2} e^{-\kappa\cos\theta} d\theta = \frac{1}{2}\int_{\theta=-\pi/4}^{\pi/4} e^{-\kappa\cos(2\theta)} d\theta = \frac{1}{2}\int_{\theta=3\pi/4}^{5\pi/4} e^{-\kappa\cos(2\theta)} d\theta = \frac{1}{2}\int_{\theta=0}^{\pi/2} e^{-\kappa\sin(2\theta)} d\theta = \int_{\theta=-\pi}^{0} e^{\kappa\sin\theta} d\theta.$$

We will need the following versions of Jordan's lemma.

### 2.3.1. *If $\mathrm{K}_A^+$ is the semicircle*

$$\mathrm{K}_A^+ = \{\lambda \in \mathbb{C} : \text{Im}\,\lambda \geq 0 \text{ and } |\lambda| = A\}, \text{ defined for } A > 0,$$

and $\mathcal{T}_A$ is a closed arc on $K_A^+$, then for every continuous function $f : \bigcup_{A>0} \mathcal{T}_A \to \mathbb{C}$, whose limit $\lim_{\lambda \to \infty} f(\lambda) = 0$, we have

$$\lim_{A \to +\infty} \int_{\mathcal{T}_A} e^{i\lambda x} f(\lambda) d\lambda = 0, \text{ when } x > 0.$$

*In particular*

$$\lim_{A \to +\infty} \int_{K_A^+} e^{i\lambda x} f(\lambda) d\lambda = 0, \text{ when } x > 0,$$

*for every continuous function* $f : \{\lambda \in \mathbb{C} : \operatorname{Im}\lambda \geq 0\} \to \mathbb{C}$ *with* $\lim_{\lambda \to \infty} f(\lambda) = 0$.

Indeed, it suffices to set $\lambda = Ae^{i\theta}$, $0 \leq \theta \leq \pi$, for $\lambda \in K_A^+$, and to notice that

$$\left| \int_{K_A^+} e^{i\lambda x} f(\lambda) d\lambda \right| = \left| \int_{\theta=0}^{\pi} e^{ixAe^{i\theta}} f(Ae^{i\theta}) Ae^{i\theta} i d\theta \right| \leq \int_{\theta=0}^{\pi} e^{-xA\sin\theta} f(Ae^{i\theta}) A d\theta \leq \frac{\pi}{x} \sup_{0 \leq \theta \leq \pi} \left| f(Ae^{i\theta}) \right|.$$

**2.3.2.** *If* $D = \{\lambda \in \mathbb{C} : \arg\lambda \in [-\pi/4, \pi/4] \cup [3\pi/4, 5\pi/4]\}$ *and* $f : D \to \mathbb{C}$ *is a continuous function with* $\lim_{\lambda \to \infty} f(\lambda) = 0$, *then*

$$\lim_{A \to +\infty} \int_{\lambda \in D \cap \{|\lambda|=A\}} \lambda e^{-\lambda^2 x} f(\lambda) d\lambda = 0, \text{ for } x > 0.$$

*More generally, if* $\mathcal{T}_A$ *are arcs on the circles* $\{\lambda \in \mathbb{C} : |\lambda| = A\}$ *with* $\mathcal{T}_A \subset D$ *(for every* $A > 0$*), then*

$$\lim_{A \to +\infty} \int_{\lambda \in \mathcal{T}_A} \lambda e^{-\lambda^2 x} f(\lambda) d\lambda = 0, \text{ for } x > 0.$$

Indeed, setting $\lambda = Ae^{i\theta}$, $\theta \in [-\pi/4, \pi/4] \cup [3\pi/4, 5\pi/4]$, for $\lambda \in D \cap \{|\lambda|=A\}$, and writing

$$\int_{\lambda \in D \cap \{|\lambda|=A\}} \lambda e^{-\lambda^2 x} f(\lambda) d\lambda = \left( \int_{\theta=-\pi/4}^{\pi/4} + \int_{\theta=3\pi/4}^{5\pi/4} \right) \left[ Ae^{i\theta} \exp\left(-A^2 x e^{i2\theta}\right) f(Ae^{i\theta}) Ai e^{i\theta} d\theta \right],$$

we find

$$\left| \int_{\lambda \in D \cap \{|\lambda|=A\}} \lambda e^{-\lambda^2 x} f(\lambda) d\lambda \right| = \left( \int_{\theta=-\pi/4}^{\pi/4} + \int_{\theta=3\pi/4}^{5\pi/4} \right) \left[ A^2 e^{-A^2 x \cos(2\theta)} \left| f(Ae^{i\theta}) \right| d\theta \right]$$

$$\leq \frac{1}{2} \left[ \sup_{\lambda \in D, |\lambda|=A} |f(\lambda)| \right] \left( \int_{\theta=-\pi/2}^{\pi/2} + \int_{\theta=3\pi/2}^{5\pi/2} \right) \left[ A^2 e^{-A^2 x \cos\theta} d\theta \right]$$

$$\leq \left[ \sup_{\lambda \in D, |\lambda|=A} |f(\lambda)| \right] \int_{\theta=0}^{\pi} A^2 e^{-A^2 x \sin\theta} d\theta \leq \frac{\pi}{x} \left[ \sup_{\lambda \in D, |\lambda|=A} |f(\lambda)| \right].$$

**2.3.3.** *Let* $m \in \mathbb{N}$ *and* $D_m = \{\lambda \in \mathbb{C} : \cos(m \arg \lambda) \geq 0\}$. *Then for a continuous function* $f : D_m \to \mathbb{C}$ *with* $\lim_{\lambda \to \infty} f(\lambda) = 0$, *we have*

$$\lim_{A \to +\infty} \int_{\lambda \in D_m \cap \{|\lambda|=A\}} \lambda^{m-1} e^{-\lambda^m x} f(\lambda) d\lambda = 0, \text{ for } x > 0.$$

The case $m=2$ is the previous one. The above more general assertion can be proved as follows: To check it first for $m=1$, and to reduce the general case to the case $m=1$, setting, in the $d\lambda$–integral, $\mu = \lambda^m$.

**2.3.4.** *Let $m \in \mathbb{N}$ and set $D_m = \{\lambda \in \mathbb{C} : \sin(m \arg \lambda) \geq 0\}$. Then for a continuous function $f : D_m \to \mathbb{C}$ with $\lim_{\lambda \to \infty} f(\lambda) = 0$, we have*

$$\lim_{A \to +\infty} \int_{\lambda \in D_m \cap \{|\lambda|=A\}} \lambda^{m-1} e^{i\lambda^m x} f(\lambda) d\lambda = 0, \text{ for } x > 0.$$

The case $m=1$ is the one in §2.3.1. The above more general assertion can be reduced to the case $m=1$, by setting, in the $d\lambda$–integral, $\mu = \lambda^m$.

**2.3.5.** *Let $D = \{\lambda \in \mathbb{C} : 0 \leq \arg \lambda \leq \pi/4\}$ and $f : D \to \mathbb{C}$ be a bounded and continuous function. Then*

$$\lim_{A \to +\infty} \int_{D \cap \{|\lambda|=A\}} \lambda^N e^{i\lambda x} e^{-\lambda^2 y} f(\lambda) d\lambda = 0, \text{ for every } N \in \mathbb{N} \cup \{0\}, \ x > 0 \text{ and } y > 0.$$

Indeed setting $\lambda = Ae^{i\theta}$, $\theta \in [0, \pi/4]$, for the points $\lambda \in D \cap \{|\lambda| = A\}$, we have

$$\left| \int_{D \cap \{|\lambda|=A\}} \lambda^N e^{i\lambda x} e^{-\lambda^2 y} f(\lambda) d\lambda \right| = \left| \int_{\theta=0}^{\pi/4} A^N e^{iN\theta} \exp(iAxe^{i\theta}) \exp(-A^2 y e^{i2\theta}) f(Ae^{i\theta}) Aie^{i\theta} d\theta \right|$$

$$\leq \int_{\theta=0}^{\pi/4} A^{N+1} e^{-Ax \sin \theta} e^{-A^2 y \cos(2\theta)} \left| f(Ae^{i\theta}) \right| d\theta$$

$$\leq \int_{\theta=0}^{\pi/8} A^{N+1} e^{-Ax \sin \theta} e^{-A^2 y \cos(2\theta)} \left| f(Ae^{i\theta}) \right| d\theta + \int_{\theta=\pi/8}^{\pi/4} A^{N+1} e^{-Ax \sin \theta} e^{-A^2 y \cos(2\theta)} \left| f(Ae^{i\theta}) \right| d\theta$$

$$\leq \left( \sup_{\lambda \in D} |f(\lambda)| \right) \left[ \int_{\theta=0}^{\pi/8} A^{N+1} e^{-Ax \sin \theta} e^{-A^2 y \cos(\pi/4)} d\theta + \int_{\theta=\pi/8}^{\pi/4} A^{N+1} e^{-Ax \sin(\pi/8)} e^{-A^2 y \cos(2\theta)} d\theta \right],$$

and the desired result follows, since

$$A^{N+1} e^{-Ax \sin \theta} e^{-A^2 y \cos(\pi/4)} \leq \frac{1}{A} \text{ for } \theta \in [0, \pi/8]$$

and

$$A^{N+1} e^{-Ax \sin(\pi/8)} e^{-A^2 y \cos(2\theta)} \leq \frac{1}{A} \text{ for } \theta \in [\pi/8, \pi/4].$$

Similarly we have

**2.3.6.** *Let $D = \{\lambda \in \mathbb{C} : 3\pi/4 \leq \arg \lambda \leq \pi\}$ and $f : D \to \mathbb{C}$ be a bounded and continuous function. Then*

$$\lim_{A \to +\infty} \int_{D \cap \{|\lambda|=A\}} \lambda^N e^{i\lambda x} e^{-\lambda^2 y} f(\lambda) d\lambda = 0, \text{ for every } N \in \mathbb{N} \cup \{0\}, \ x > 0 \text{ and } y > 0.$$

**2.3.7.** *Let $D = \{\lambda \in \mathbb{C} : \pi/4 \leq \arg \lambda \leq 3\pi/4\}$ and $f : D \to \mathbb{C}$ be a bounded and continuous function. Then*

$$\lim_{A \to +\infty} \int_{D \cap \{|\lambda|=A\}} \lambda^N e^{i\lambda x} e^{\lambda^2 y} f(\lambda) d\lambda = 0, \text{ for every } N \in \mathbb{N} \cup \{0\}, \ x > 0 \text{ and } y > 0.$$

**2.4. Elementary integration by parts formulas.** A. Let $u_0(x) \in \mathcal{S}([0,\infty))$. Then for $\lambda \in \mathbb{C}$ with $\operatorname{Im}\lambda \leq 0$ and $\lambda \neq 0$, we have

$$\hat{u}_0(\lambda) = \frac{u_0(0)}{i\lambda} + \frac{1}{i\lambda} \int_{y=0}^{\infty} e^{-i\lambda x} u_0'(x) dx$$

$$= \frac{u_0(0)}{i\lambda} + \frac{u_0'(0)}{(i\lambda)^2} + \frac{1}{(i\lambda)^2} \int_{y=0}^{\infty} e^{-i\lambda x} u_0''(x) dx$$

$$= \frac{u_0(0)}{i\lambda} + \frac{u_0'(0)}{(i\lambda)^2} + \frac{u_0''(0)}{(i\lambda)^3} + \frac{1}{(i\lambda)^3} \int_{y=0}^{\infty} e^{-i\lambda x} u_0'''(x) dx = \sum_{n=0}^{m} \frac{u_0^{(n)}(0)}{(i\lambda)^n} + \frac{1}{(i\lambda)^m} \int_{y=0}^{\infty} e^{-i\lambda x} u_0^{(m+1)}(x) dx. \quad (2.4)$$

In particular,

$$|\hat{u}_0(\lambda)| \leq \frac{|u_0(0)|}{|\lambda|} + \frac{1}{|\lambda|} \int_{y=0}^{\infty} |u_0'(x)| dx = O\left(\frac{1}{\lambda}\right) \text{ as } \lambda \to \infty \text{ with } \lambda \in \mathbb{C} \text{ and } \operatorname{Im}\lambda \leq 0. \quad (2.5)$$

Also

$$\frac{1}{(i\lambda)^m} \int_{y=0}^{\infty} e^{-i\lambda x} u_0^{(m+1)}(x) dx = \frac{u_0^{(m+1)}(0)}{(i\lambda)^{m+1}} + \frac{1}{(i\lambda)^{m+1}} \int_{y=0}^{\infty} e^{-i\lambda x} u_0^{(m+2)}(x) dx,$$

and therefore, the last term in (2.1) is estimated by

$$\left| \frac{1}{(i\lambda)^m} \int_{y=0}^{\infty} e^{-i\lambda x} u_0^{(m+1)}(x) dx \right| = O\left(\frac{1}{\lambda^{m+1}}\right) \text{ as } \lambda \to \infty \text{ with } \lambda \in \mathbb{C} \text{ and } \operatorname{Im}\lambda \leq 0. \quad (2.6)$$

B. Let $g_0(t) \in C^{\infty}([0,\infty))$. For $\lambda \in \mathbb{C}$, $\lambda \neq 0$, we have

$$e^{-\lambda^2 t} \int_{\tau=0}^{t} e^{\lambda^2 \tau} g_0(\tau) d\tau = \frac{g_0(t)}{\lambda^2} - \frac{g_0(0)}{\lambda^2} e^{-\lambda^2 t} - \frac{1}{\lambda^2} e^{-\lambda^2 t} \int_{\tau=0}^{t} e^{\lambda^2 \tau} g_0'(\tau) d\tau$$

$$= \frac{g_0(t)}{\lambda^2} - \frac{g_0(0)}{\lambda^2} e^{-\lambda^2 t} - \frac{g_0'(t)}{\lambda^4} + \frac{g_0'(0)}{\lambda^4} e^{-\lambda^2 t} + \frac{1}{\lambda^4} e^{-\lambda^2 t} \int_{\tau=0}^{t} e^{\lambda^2 \tau} g_0''(\tau) d\tau$$

$$= \sum_{n=0}^{m} (-1)^n \left[ \frac{g_0^{(n)}(t)}{\lambda^{2n+2}} - \frac{g_0^{(n)}(0)}{\lambda^{2n+2}} e^{-\lambda^2 t} \right] + (-1)^{m+1} \frac{1}{\lambda^{2m+2}} e^{-\lambda^2 t} \int_{\tau=0}^{t} e^{\lambda^2 \tau} g_0^{(m+1)}(\tau) d\tau. \quad (2.7)$$

It follows that

$$\left| e^{-\lambda^2 t} \int_{\tau=0}^{t} e^{\lambda^2 \tau} g_0(\tau) d\tau \right| \leq \frac{|g_0(t)|}{|\lambda|^2} + \frac{|g_0(0)|}{|\lambda|^2} \left| e^{-\lambda^2 t} \right| + \frac{1}{|\lambda|^2} \left| e^{-\lambda^2 t} \int_{\tau=0}^{t} e^{\lambda^2 \tau} g_0'(\tau) d\tau \right|$$

$$\leq \frac{|g_0(t)|}{|\lambda|^2} + \frac{|g_0(0)|}{|\lambda|^2} \left| e^{-\lambda^2 t} \right| + \frac{1}{|\lambda|^2} \int_{\tau=0}^{t} \left| e^{-\lambda^2 (t-\tau)} g_0'(\tau) \right| d\tau. \quad (2.8)$$

In particular, if $\operatorname{Re}(\lambda^2) \geq 0$ then $\left| e^{-\lambda^2(t-\tau)} \right| = e^{-(\operatorname{Re}\lambda^2)(t-\tau)} \leq 1$ for $t \geq \tau$, and therefore, for a fixed $t > 0$,

$$\left| e^{-\lambda^2 t} \int_{\tau=0}^{t} e^{\lambda^2 \tau} g_0(\tau) d\tau \right|$$

$$\leq \frac{|g_0(t)|}{|\lambda|^2} + \frac{|g_0(0)|}{|\lambda|^2} + \frac{1}{|\lambda|^2} \sup_{0 \leq \tau \leq t} |g_0'(\tau)| = O\left(\frac{1}{\lambda^2}\right) \text{ as } \lambda \to \infty \text{ with } \lambda \in \mathbb{C} \text{ and } \operatorname{Re}(\lambda^2) \geq 0. \quad (2.9)$$

Applying (2.6) with $g_0^{(m)}(t)$ in place of $g_0(t)$, we have the following estimate for the last term in (2.4):

$$\left| \frac{1}{\lambda^{2m}} e^{-\lambda^2 t} \int_{\tau=0}^{t} e^{\lambda^2 \tau} g_0^{(m)}(\tau) d\tau \right| = O\left(\frac{1}{\lambda^{2m+2}}\right) \text{ as } \lambda \to \infty \text{ with } \lambda \in \mathbb{C} \text{ and } \operatorname{Re}(\lambda^2) \geq 0. \quad (2.10)$$

## 3. Derivation of the formula

In this section we outline the derivation of the formula (1.2) by the method of Fokas. In searching for the solution $u(x,t)$, we will assume that it exists and, moreover, is sufficiently smooth and appropriately rapidly decreasing as $x \to +\infty$. We will also make some further assumptions, which we will point out in the process of derivation.

***Step 1*** Setting

$$\hat{u}(\lambda,t) = \int_{x=0}^{\infty} u(x,t)e^{-i\lambda x}dx,$$

we see that $\hat{u}(\lambda,t)$ is an analytic function for $\lambda \in \mathbb{C}$ with $\operatorname{Im}\lambda < 0$ (for each fixed $t > 0$) with continuous extension upto the line $\{\lambda \in \mathbb{C}: \operatorname{Re}\lambda = 0\}$. This follows from Lebesque's dominated convergence theorem (case (2.3)) and the Cauchy-Riemann equations, using the condition $\sup_{x\geq 0}|x^3 u(x,t)| < \infty$. Then, again by (2.3),

$$\frac{\partial \hat{u}(\lambda,t)}{\partial t} = \int_{x=0}^{\infty} e^{-i\lambda x} \frac{\partial u(x,t)}{\partial t} dx \ (t > 0),$$

using the condition that $\sup_{x\geq 0, \alpha \leq t \leq \beta}\left|x^2 \frac{\partial u(x,t)}{\partial x}\right| < \infty$ (for $\beta > \alpha > 0$), which guarantees that

$$\sup_{\alpha \leq t \leq \beta} \int_{x=0}^{\infty} \left|\frac{\partial u(x,t)}{\partial t}\right| dx < \infty.$$

Thus the equation $u_t = u_{xx}$ gives

$$\frac{\partial \hat{u}(\lambda,t)}{\partial t} = \int_{x=0}^{\infty} e^{-i\lambda x} \frac{\partial^2 u(x,t)}{\partial x^2} dx \ (\lambda \in \mathbb{C} \text{ with } \operatorname{Im}\lambda \leq 0).$$

Integration by parts leads to the differential equation

$$\hat{u}_t + \lambda^2 \hat{u} = -g_1(t) - i\lambda g_0(t) \ (t > 0),$$

where $\hat{u}_t = \frac{\partial \hat{u}(\lambda,t)}{\partial t}$ και $g_1(t) = \frac{\partial u(0,t)}{\partial x}$. Solving the above equation we find that

$$\hat{u}(\lambda,t) = e^{-\lambda^2 t}\hat{u}_0(\lambda) - e^{-\lambda^2 t}[\tilde{g}_1(\lambda,t) + i\lambda \tilde{g}_0(\lambda,t)], \ \lambda \in \mathbb{C} \text{ with } \operatorname{Im}\lambda \leq 0, \quad (3.1)$$

where $\tilde{g}_1(\lambda,t) = \int_{\tau=0}^{t} e^{\lambda^2 \tau} g_1(\tau)d\tau$.

(Comment 1: At this point we need to have $g_1 \in C([0,\infty))$, i.e., continuous upto the point $t = 0$, and that

$$\hat{u}(\lambda,t) = \int_{x=0}^{\infty} e^{-i\lambda x} u(x,t)dx \to \int_{x=0}^{\infty} e^{-i\lambda x} u(x,0)dx = \hat{u}_0(\lambda), \text{ as } t \to 0^+.)$$

By Fourier inversion formula (2.2),

$$u(x,t) = \frac{1}{2\pi} \int_{\lambda=-\infty}^{\infty} \hat{u}(\lambda,t)e^{i\lambda x}d\lambda$$

$$= \frac{1}{2\pi} \int_{\lambda=-\infty}^{\infty} e^{i\lambda x - \lambda^2 t}\hat{u}_0(\lambda)d\lambda - \frac{1}{2\pi} \int_{\lambda=-\infty}^{\infty} e^{i\lambda x - \lambda^2 t}[\tilde{g}_1(\lambda,t) + i\lambda \tilde{g}_0(\lambda,t)]d\lambda, \text{ for } x > 0, t > 0. \quad (3.2)$$

(Comment 2: In the above equations, the first integral is interpreted as the limit $\lim_{A \to \infty} \int_{-A}^{A}$, and so is the third one. Observe that the second integral – moreover – converges absolutely because of the factor $e^{-\lambda^2 t}$.)

**Step 2** Integrating by parts, as in (2.4), we obtain

$$\widetilde{g}_1(\lambda,t) = \int_{\tau=0}^{t} e^{\lambda^2\tau} g_1(\tau) d\tau = \frac{1}{\lambda^2} e^{\lambda^2 t} g_1(t) - \frac{1}{\lambda^2} g_1(0) - \frac{1}{\lambda^2} \int_{\tau=0}^{t} e^{\lambda^2\tau} \frac{dg_1(\tau)}{d\tau} d\tau.$$

(Comment 3: At this point of the derivation process we need to assume that $g_1 \in C^1([0,\infty))$.)
It follows that, for each fixed $t > 0$,

$$G_1(\lambda,t) \stackrel{def}{=} e^{-\lambda^2 t} \widetilde{g}_1(\lambda,t) = O\left(\frac{1}{\lambda^2}\right) \text{ as } \lambda \to \infty \text{ with } \lambda \in \Omega^+, \tag{3.3}$$

where $\Omega^+ \stackrel{def}{=} \{\lambda \in \mathbb{C}: \text{Im}\,\lambda \geq 0 \text{ and } \text{Re}(\lambda^2) \geq 0\}$. Notice that if $\lambda = \xi + i\eta$ with $\xi,\eta \in \mathbb{R}$, then

$$\lambda^2 = (\xi^2 - \eta^2) + 2\xi\eta i, \quad \Omega^+ = \{\xi + i\eta : \eta \geq 0 \text{ and } \xi^2 - \eta^2 \geq 0\},$$

and

$$\text{for } \lambda \in \Omega^+, \left|e^{-\lambda^2 t}\right| = e^{-\text{Re}(\lambda^2)t} \leq 1 \text{ and } \left|e^{-\lambda^2(t-\tau)}\right| \leq 1 \text{ for } \tau \leq t.$$

Similarly, for each fixed $t > 0$,

$$G_0(\lambda,t) \stackrel{def}{=} e^{-\lambda^2 t} \widetilde{g}_0(\lambda,t) = O\left(\frac{1}{\lambda^2}\right) \text{ as } \lambda \to \infty \text{ with } \lambda \in \Omega^+. \tag{3.4}$$

We claim that

$$\int_{\lambda=-\infty}^{\infty} e^{i\lambda x}[G_1(\lambda,t) + i\lambda G_0(\lambda,t)]d\lambda = \int_{\lambda \in \Gamma} e^{i\lambda x}[G_1(\lambda,t) + i\lambda G_0(\lambda,t)]d\lambda \quad (x > 0). \tag{3.5}$$

Indeed, since the function $e^{i\lambda x}[G_1(\lambda,t) + i\lambda G_0(\lambda,t)]$ is analytic for $\lambda \in \mathbb{C}$, by Cauchy's theorem:

$$\int_{\lambda=-A}^{A} e^{i\lambda x}[G_1(\lambda,t) + i\lambda G_0(\lambda,t)]d\lambda - \int_{\lambda \in \Gamma_A} e^{i\lambda x}[G_1(\lambda,t) + i\lambda G_0(\lambda,t)]d\lambda$$

$$+ \int_{\lambda \in K(0,A) \cap \Omega^+} e^{i\lambda x}[G_1(\lambda,t) + i\lambda G_0(\lambda,t)]d\lambda = 0 \tag{3.6}$$

where $K(0,A) \cap \Omega^+$ is the part of the circle $K(0,A) = \{\lambda \in \mathbb{C}: |\lambda| \leq A\}$, $A > 0$, which lies inside the set $\Omega^+$ and $\Gamma_A$ is the part of $\Gamma$ which lies inside this circle.

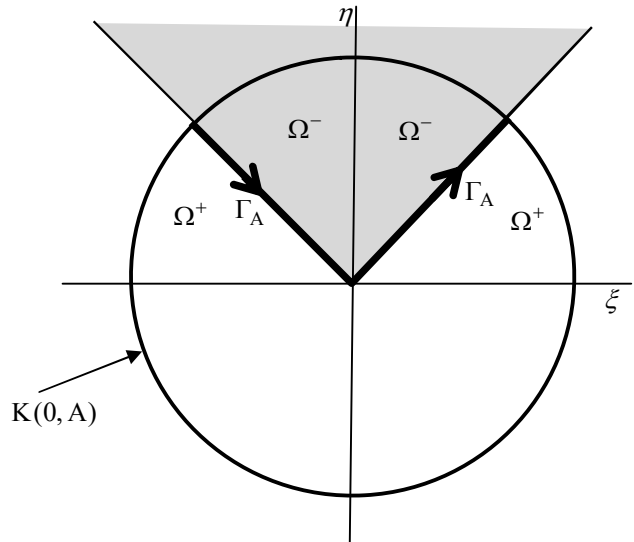

Now, by (3.3) and (3.4),

$$G_1(\lambda,t) + i\lambda G_0(\lambda,t) = O\left(\frac{1}{\lambda}\right) \text{ as } \lambda \to \infty \text{ with } \lambda \in \Omega^+,$$

whence, from Jordan's lemma 2.3.1,

Andreas Chatziafratis$$\lim_{A \to \infty} \int_{\lambda \in K(0,A) \cap \Omega^+} e^{i\lambda x}[G_1(\lambda,t) + i\lambda G_0(\lambda,t)]d\lambda = 0.$$

Therefore (3.5) follows from (3.6), if we let $A \to \infty$.

*Remark* In the proof of (3.5), the integral $\int_{\lambda \in \Gamma}$ is interpreted as $\lim_{A \to \infty} \int_{\lambda \in \Gamma_A}$. However, because of the factor $e^{i\lambda x}$, the integral $\int_{\lambda \in \Gamma}$ converges absolutely, whereas the integral $\int_{-\infty}^{\infty}$ in (3.5), in general, does not converge absolutely. (The absolute convergence of $\int_{\lambda \in \Gamma}$ follows from the fact that, for $\lambda \in \Gamma$, $\left|e^{i\lambda x}\right| = e^{-x \operatorname{Im} \lambda} = e^{-x|\lambda|/\sqrt{2}}$.)

***Step 3*** In view of (3.5), (3.2) can be written as follows:

$$u(x,t) = \frac{1}{2\pi} \int_{\lambda=-\infty}^{\infty} e^{i\lambda x - \lambda^2 t} \hat{u}_0(\lambda) d\lambda - \frac{1}{2\pi} \int_{\lambda \in \Gamma} e^{i\lambda x - \lambda^2 t}[\tilde{g}_1(\lambda,t) + i\lambda \tilde{g}_0(\lambda,t)]d\lambda, \text{ for } x > 0, \ t > 0. \quad (3.7)$$

Next setting $-\lambda$ in place of $\lambda$ in equation (3.1), we obtain

$$\hat{u}(-\lambda,t) = e^{-\lambda^2 t} \hat{u}_0(-\lambda) - e^{-\lambda^2 t}[\tilde{g}_1(\lambda,t) - i\lambda \tilde{g}_0(\lambda,t)], \text{ for } \lambda \in \mathbb{C} \text{ with } \operatorname{Im} \lambda \geq 0. \quad (3.8)$$

Integrating (11), we have

$$\int_{\lambda \in \Gamma} e^{i\lambda x - \lambda^2 t} \tilde{g}_1(\lambda,t) d\lambda = -\int_{\lambda \in \Gamma} e^{i\lambda x} \hat{u}(-\lambda,t) d\lambda + \int_{\lambda \in \Gamma} e^{i\lambda x - \lambda^2 t} \hat{u}_0(-\lambda) d\lambda + \int_{\lambda \in \Gamma} e^{i\lambda x - \lambda^2 t} i\lambda \tilde{g}_0(\lambda,t) d\lambda. \quad (3.9)$$

We claim that

$$\int_{\lambda \in \Gamma} e^{i\lambda x} \hat{u}(-\lambda,t) d\lambda = 0. \quad (3.10)$$

Indeed, since $e^{i\lambda x} \hat{u}(-\lambda,t)$ is continuous in $\{\lambda \in \mathbb{C} : \operatorname{Im} \lambda \geq 0\}$ and analytic in $\{\lambda \in \mathbb{C} : \operatorname{Im} \lambda > 0\}$, by Cauchy's theorem

$$\int_{\lambda \in \Gamma_A} e^{i\lambda x} \hat{u}(-\lambda,t) d\lambda + \int_{\lambda \in K(0,A) \cap \Omega^-} e^{i\lambda x} \hat{u}(-\lambda,t) d\lambda = 0 \text{ for } A > 0. \quad (3.11)$$

On the other hand, integration by parts, as in (2.1), gives

$$\hat{u}(-\lambda,t) = \int_{x=0}^{\infty} u(x,t) e^{i\lambda x} dx = -\frac{1}{i\lambda} u(0,t) - \int_{x=0}^{\infty} \frac{1}{i\lambda} \frac{\partial u(x,t)}{\partial x} e^{i\lambda x} dx$$

which implies that

$$\hat{u}(-\lambda,t) = O\left(\frac{1}{\lambda}\right) \text{ as } \lambda \to \infty \text{ for } \lambda \in \mathbb{C} \text{ with } \operatorname{Im} \lambda \geq 0.$$

($t > 0$ is kept fixed.) Therefore letting $A \to \infty$ and using Jordan's lemma 2.3.1, (3.11) gives (3.10). Finally, substituting (3.10) and (3.9) in (3.7), we obtain (1.2). □

## 4. Other forms of the solution

### 4.1. Writing the solution with an Ehrenpreis integral for $0 < t < T$.

Let us fix a $T > 0$. We will show that for $0 < t < T$, the solution given by (1.2) can be written – equivalently – as follows: For $0 < t < T$ and $x > 0$,

$$u(x,t) = \frac{1}{2\pi} \int_{\lambda=-\infty}^{\infty} e^{i\lambda x - \lambda^2 t} \hat{u}_0(\lambda) d\lambda - \frac{1}{2\pi} \int_{\lambda \in \Gamma} e^{i\lambda x - \lambda^2 t} \hat{u}_0(-\lambda) d\lambda - \frac{i}{\pi} \int_{\lambda \in \Gamma} e^{i\lambda x - \lambda^2 t} \lambda \tilde{g}_0(\lambda, T) d\lambda. \quad (4.1)$$

Indeed, the difference

$$\delta(x,t) = \frac{i}{\pi} \int_{\lambda \in \Gamma} e^{i\lambda x - \lambda^2 t} \lambda [\widetilde{g}_0(\lambda, T) - \widetilde{g}_0(\lambda, t)] d\lambda = \frac{i}{\pi} \int_{\lambda \in \Gamma} e^{i\lambda x - \lambda^2 t} \lambda \left[ \int_{\tau=t}^{T} e^{\lambda^2 \tau} g_0(\tau) d\tau \right] d\lambda$$

$$= \frac{i}{\pi} \int_{\lambda \in \Gamma} e^{i\lambda x} \lambda \left[ \int_{\tau=t}^{T} e^{\lambda^2 (\tau-t)} g_0(\tau) d\tau \right] d\lambda.$$

But, integration by parts gives

$$\int_{\tau=t}^{T} e^{\lambda^2(\tau-t)} g_0(\tau) d\tau = \frac{1}{\lambda^2} \int_{\tau=t}^{T} \frac{d[e^{\lambda^2(\tau-t)}]}{d\tau} g_0(\tau) d\tau = \frac{1}{\lambda^2} \left[ e^{\lambda^2(\tau-t)} g_0(\tau) \right]_{\tau=t}^{T} - \frac{1}{\lambda^2} \int_{\tau=t}^{T} e^{\lambda^2(\tau-t)} \frac{dg_0(\tau)}{d\tau} d\tau$$

$$= \frac{1}{\lambda^2} e^{\lambda^2(T-t)} g_0(T) - \frac{1}{\lambda^2} g_0(t) - \frac{1}{\lambda^2} \int_{\tau=t}^{T} e^{\lambda^2(\tau-t)} \frac{dg_0(\tau)}{d\tau} d\tau.$$

Since for $\lambda \in \Omega^-$ and $\tau \geq t$, $e^{\operatorname{Re}(\lambda^2)(\tau-t)} \leq 1$, it follows that

$$\lambda \int_{\tau=t}^{T} e^{\lambda^2(\tau-t)} g_0(\tau) d\tau = O\left(\frac{1}{\lambda}\right) \text{ as } \lambda \to \infty \text{ with } \lambda \in \Omega^- \text{ (} t \text{ is kept fixed, } 0 < t < T \text{)}.$$

Therefore by Cauchy's theorem in $\Omega^-$ and Jordan's lemma 2.3.1,

$$\int_{\lambda \in \Gamma} e^{i\lambda x} \lambda \left[ \int_{\tau=t}^{T} e^{\lambda^2(\tau-t)} g_0(\tau) d\tau \right] d\lambda = \lim_{A \to \infty} \int_{\Gamma_A} \cdots = 0,$$

i.e., $\delta(x,t) = 0$. This completes the proof of (4.1).

Thus (4.1) gives the solution $u(x,t)$ as a combination – integral of the exponential solutions $e^{i\lambda x - \lambda^2 t}$, $\lambda \in \mathbb{C}$, of the equation, with measures being independent of $x$ and $t$.

### 4.2. Writing the solution with an integral involving Gauss's kernel.

*Proposition 1* For a continuous function $g_0(t)$, defined for $t \geq 0$, with $\int_0^{\infty} |g_0(t)| dt < \infty$, we have

$$-\frac{i}{\pi} \int_{\lambda \in \Gamma} e^{i\lambda x - \lambda^2 t} \lambda \widetilde{g}_0(\lambda, t) d\lambda = \frac{1}{2\sqrt{\pi}} \int_{\tau=0}^{t} g_0(\tau) \frac{x}{(t-\tau)^{3/2}} \exp\left(-\frac{x^2}{4(t-\tau)}\right) d\tau \quad (x > 0, \ t > 0).$$

*Proof* Since

$$\left| e^{i\lambda x - \lambda^2 t} \lambda e^{\lambda^2 \tau} \right| = e^{-x \operatorname{Im} \lambda} |\lambda| = e^{-x|\lambda|/\sqrt{2}} |\lambda| \quad \text{when } \lambda \in \Gamma,$$

the double integral, corresponding to the iterated integral

$$\int_{\lambda \in \Gamma} \int_{\tau=0}^{t} e^{i\lambda x - \lambda^2 t} \lambda e^{\lambda^2 \tau} g_0(\tau) d\tau d\lambda,$$

converges absolutely (for $x > 0$, $t > 0$). Therefore, by Fubini's theorem,

$$\int_{\lambda \in \Gamma} e^{i\lambda x - \lambda^2 t} \lambda \widetilde{g}_0(\lambda, t) d\lambda = \int_{\tau=0}^{t} g_0(\tau) \left( \int_{\lambda \in \Gamma} e^{i\lambda x} e^{-\lambda^2(t-\tau)} \lambda d\lambda \right) d\tau.$$

But, from Cauchy's theorem and Jordan's lemmas 2.3.5 and 2.3.6,

$$\int_{\lambda \in \Gamma} e^{i\lambda x} e^{-\lambda^2(t-\tau)} \lambda d\lambda = \int_{\lambda=-\infty}^{\infty} e^{i\lambda x} e^{-\lambda^2(t-\tau)} \lambda d\lambda \quad \text{for } \tau < t.$$

Also, recalling the Fourier transform of the function $e^{-\kappa \lambda^2}$ (considered as a function of $\lambda$), where $\kappa$ is a positive constant, i.e.,

$$\int_{\lambda=-\infty}^{\infty} e^{-i\lambda x} e^{-\kappa \lambda^2} d\lambda = \sqrt{\frac{\pi}{\kappa}} e^{-x^2/4\kappa},$$

and differentiating with respect to $x$, we find

$$\int_{\lambda=-\infty}^{\infty} e^{-i\lambda x} e^{-\kappa \lambda^2} \lambda \, d\lambda = i \frac{\sqrt{\pi}}{2} \frac{x}{\kappa^{3/2}} e^{-x^2/4\kappa} \quad (\kappa > 0).$$

Therefore

$$\int_{\lambda=-\infty}^{\infty} e^{i\lambda x} e^{-\lambda^2(t-\tau)} \lambda \, d\lambda = -i \frac{\sqrt{\pi}}{2} \frac{x}{(t-\tau)^{3/2}} \exp\left(\frac{-x^2}{4(t-\tau)}\right), \text{ for } \tau < t,$$

and the formula of the proposition follows. □

*Proposition 2* For a continuous function $u_0(s)$, defined for $s \geq 0$, with $\int_0^\infty |u_0(s)| ds < \infty$, we have

$$\frac{1}{2\pi} \left( \int_{\lambda=-\infty}^{\infty} e^{i\lambda x - \lambda^2 t} \hat{u}_0(\lambda) d\lambda - \int_{\lambda \in \Gamma} e^{i\lambda x - \lambda^2 t} \hat{u}_0(-\lambda) d\lambda \right) = \frac{1}{2\sqrt{\pi t}} \int_{s=0}^{\infty} u_0(s) \left[ \exp\left[-\frac{(s-x)^2}{4t}\right] - \exp\left[-\frac{(s+x)^2}{4t}\right] \right] ds.$$

*Proof.* It follows from Fubini's theorem that, for $x > 0$ and $t > 0$,

$$\int_{\lambda=-\infty}^{\infty} e^{i\lambda x - \lambda^2 t} \hat{u}_0(\lambda) d\lambda = \int_{\lambda=-\infty}^{\infty} e^{i\lambda x - \lambda^2 t} \left( \int_{s=0}^{\infty} e^{-i\lambda s} u_0(s) ds \right) d\lambda = \int_{s=0}^{\infty} u_0(s) \left( \int_{\lambda=-\infty}^{\infty} e^{-i\lambda(s-x)} e^{-\lambda^2 t} d\lambda \right) ds,$$

and, in combination with Cauchy's theorem and Jordan's lemmas 2.3.5 and 2.3.6,

$$\int_{\lambda \in \Gamma} e^{i\lambda x - \lambda^2 t} \hat{u}_0(-\lambda) d\lambda = \int_{\lambda \in \Gamma} e^{i\lambda x - \lambda^2 t} \left( \int_{s=0}^{\infty} e^{i\lambda s} u_0(s) ds \right) d\lambda = \int_{s=0}^{\infty} u_0(s) \left( \int_{\lambda \in \Gamma} e^{i\lambda s} e^{i\lambda x - \lambda^2 t} d\lambda \right) ds$$

$$= -\int_{s=0}^{\infty} u_0(s) \left( \int_{\lambda=-\infty}^{\infty} e^{-i\lambda(s+x)} e^{-\lambda^2 t} d\lambda \right) ds.$$

Thus, the equation of the proposition follows from the formula of the Fourier transform of the function $e^{-\kappa \lambda^2}$ (of $\lambda$). □

*Conclusion* The function $u(x,t)$ which is defined by (1.2) – equivalently by (4.1) – for $x > 0$ and $t > 0$, can be written in the form

$$u(x,t) = \frac{1}{2\sqrt{\pi t}} \int_{s=0}^{\infty} u_0(s) \left[ \exp\left[-\frac{(s-x)^2}{4t}\right] - \exp\left[-\frac{(s+x)^2}{4t}\right] \right] ds$$

$$+ \frac{1}{2\sqrt{\pi}} \int_{\tau=0}^{t} g_0(\tau) \frac{x}{(t-\tau)^{3/2}} \exp\left(-\frac{x^2}{4(t-\tau)}\right) d\tau. \quad (4.2)$$

*Remark* It is immediately verified that the above function $u(x,t)$ – in view of Lebesque's dominated convergence theorem – satisfies the equation $u_t = u_{xx}$, for $x > 0$ and $t > 0$. Let us notice however that the second integral in (4.2) for $x = 0$ is equal to zero, while its limit as $x \to 0^+$, as we will show, is equal to $g_0(t)$. There is a difficulty also with the limit of the first integral in (4.2), as $t \to 0^+$, which, as we will see is equal to $u_0(x)$. (See [1] for the first integral and [14, Theorem 62.2] for the second integral of this form of the solution. See also [16] for a detailed study of the heat equation.)

### 4.3. Writing the solution with integrals taken over the real line.

Equation (1.2) can be written as follows: For $x > 0$ and $t > 0$,

$$u(x,t) = \frac{1}{2\pi} \int_{\lambda=-\infty}^{\infty} e^{i\lambda x - \lambda^2 t} \hat{u}_0(\lambda) d\lambda - \frac{1}{2\pi} \int_{\lambda=-\infty}^{\infty} e^{i\lambda x - \lambda^2 t} \hat{u}_0(-\lambda) d\lambda - \frac{i}{\pi} \int_{\lambda=-\infty}^{\infty} e^{i\lambda x - \lambda^2 t} \lambda \tilde{g}_0(\lambda, t) d\lambda. \quad (4.3)$$

Observe that while the first two integrals in the RHS of the above equation converge absolutely, the third one, i.e., the integral

$$\int_{\lambda=-\infty}^{\infty} e^{i\lambda x - \lambda^2 t} \lambda \widetilde{g}_0(\lambda, t) d\lambda = \lim_{A \to \infty} \int_{\lambda=-A}^{A} e^{i\lambda x - \lambda^2 t} \lambda \widetilde{g}_0(\lambda, t) d\lambda,$$

does not converge absolutely – in general – and it is interpreted as the above limit. Notice also

$$\left[ \int_{\lambda=-\infty}^{\infty} e^{i\lambda x - \lambda^2 t} \lambda \widetilde{g}_0(\lambda, t) d\lambda \right]_{x=0} = 0.$$

Indeed, for every $A > 0$,

$$\left[ \int_{\lambda=-A}^{A} e^{i\lambda x - \lambda^2 t} \lambda \widetilde{g}_0(\lambda, t) d\lambda \right]_{x=0} = \int_{\lambda=-A}^{A} e^{-\lambda^2 t} \lambda \left( \int_{\tau=0}^{t} e^{-\lambda^2 \tau} g_0(\tau) d\tau \right) d\lambda = 0,$$

since the function $e^{-\lambda^2 t} \lambda \left( \int_{\tau=0}^{t} e^{-\lambda^2 \tau} g_0(\tau) d\tau \right)$ is odd – with respect to $\lambda \in \mathbb{R}$.

But, as we will show,

$$\lim_{x \to 0^+} \int_{\lambda=-\infty}^{\infty} e^{i\lambda x - \lambda^2 t} \lambda \widetilde{g}_0(\lambda, t) d\lambda = \lim_{x \to 0^+} \left( \lim_{A \to \infty} \int_{\lambda=-A}^{A} e^{i\lambda x - \lambda^2 t} \lambda \widetilde{g}_0(\lambda, t) d\lambda \right) = -\frac{\pi}{i} g_0(t).$$

Finally it is easy to check that (4.3) can be written also in the following way:

$$u(x,t) = \frac{2}{\pi} \int_{\lambda=0}^{\infty} \sin(\lambda x) e^{-\lambda^2 t} \left[ \int_{y=0}^{\infty} \sin(\lambda y) u_0(y) dy + \lambda \int_{\tau=0}^{t} e^{\lambda^2 \tau} g_0(\tau) d\tau \right] d\lambda \quad (x > 0, t > 0). \quad (4.4).$$

## 5. Proof of Theorem 1

We split the proof in several steps.

***Step 1*** We claim that the integrals in the RHS of (4.1) are absolutely and uniformly convergent on compact subsets of $Q$ and remain so after any number of differentiations – with respect to $x$ or $t$.
Firstly, the convergence

$$\int_{\lambda=-\infty}^{\infty} \left| \lambda^N e^{i\lambda x - \lambda^2 t} \hat{u}_0(\lambda) \right| d\lambda < \infty \quad (x > 0, t > 0 \text{ and } N \in \mathbb{N}) \quad (5.1)$$

follows from the presence of the factor $e^{-\lambda^2 t}$ (since $t > 0$) and the fact that the function $\hat{u}_0(\lambda)$ is bounded for $\lambda \in \mathbb{R}$.
Also the integral

$$\int_{\lambda \in \Gamma} \lambda^N e^{i\lambda x - \lambda^2 t} \hat{u}_0(-\lambda) d\lambda \quad (x > 0, t > 0) \quad (5.2)$$

converges absolutely, for every $N \in \mathbb{N}$. This follows from the fact that

$$\left| e^{i\lambda x - \lambda^2 t} \right| = e^{-x \operatorname{Im} \lambda} = e^{-x|\lambda|/\sqrt{2}} \quad \text{for } \lambda \in \Gamma,$$

and

$$\left| \hat{u}_0(-\lambda) \right| = \left| \int_{x=0}^{\infty} u_0(x) e^{i\lambda x} dx \right| \leq \int_{x=0}^{\infty} |u_0(x)| dx \quad (\text{since } \left| e^{i\lambda x} \right| \leq 1) \text{ for } \lambda \in \Gamma.$$

Similarly the integral

$$\int_{\lambda \in \Gamma} \lambda^N e^{i\lambda x - \lambda^2 t} \lambda \widetilde{g}_0(\lambda, T) d\lambda \quad (0 < t < T) \quad (5.3)$$

Andreas Chatziafratisconverges absolutely for $x > 0$ and for every $N \in \mathbb{N}$. This again follows from the presence of the factor $e^{i\lambda x}$ and the fact that for $\lambda \in \Gamma$, $\left|e^{-\lambda^2 t}\right| = \left|e^{-\lambda^2 \tau}\right| = 1$, so that the function $e^{-\lambda^2 t}\tilde{g}_0(\lambda, T)$ remains bounded for $\lambda \in \Gamma$:

$$\left|e^{-\lambda^2 t}\tilde{g}_0(\lambda, T)\right| = \left|e^{-\lambda^2 t}\int_{\tau=0}^{T} e^{\lambda^2 \tau} g_0(\tau) d\tau\right| \leq \int_{\tau=0}^{T} |g_0(\tau)| d\tau.$$

In addition, it is a easy to check that, given $x_0 > 0$ and $t_0 > 0$, the absolute convergence of the integrals in (5.1), (5.2) and (5.3), is uniform for $x \geq x_0$ and $t \geq t_0$, and this implies the claim.

**Step 2** For $x > 0$ and $t > 0$,

$$\frac{\partial}{\partial t}\left(\int_{\lambda=-\infty}^{\infty} e^{i\lambda x - \lambda^2 t} \hat{u}_0(\lambda) d\lambda\right) = \int_{\lambda=-\infty}^{\infty} \frac{\partial [e^{i\lambda x - \lambda^2 t}]}{\partial t} \hat{u}_0(\lambda) d\lambda.$$

Indeed, this follows from Lebesgue's dominated convergence theorem (2.3), since

$$\sup_{t \geq t_0} \int_{\lambda=-\infty}^{\infty} \left|\lambda^2 e^{i\lambda x - \lambda^2 t} \hat{u}_0(\lambda)\right| d\lambda < \infty, \text{ for } t_0 > 0.$$

Similarly, for $x > 0$ and $t > 0$,

$$\frac{\partial^2}{\partial x^2}\left(\int_{\lambda=-\infty}^{\infty} e^{i\lambda x - \lambda^2 t} \hat{u}_0(\lambda) d\lambda\right) = \int_{\lambda=-\infty}^{\infty} \frac{\partial^2 [e^{i\lambda x - \lambda^2 t}]}{\partial x^2} \hat{u}_0(\lambda) d\lambda.$$

At this point we use the fact that

$$\sup_{x > 0} \int_{\lambda=-\infty}^{\infty} \left|\lambda^2 e^{i\lambda x - \lambda^2 t} \hat{u}_0(\lambda)\right| d\lambda < \infty, \text{ for } t > 0.$$

In general, for $x > 0$ and $t > 0$,

$$\frac{\partial^{k+l}}{\partial t^k \partial x^l}\left(\int_{\lambda=-\infty}^{\infty} e^{i\lambda x - \lambda^2 t} \hat{u}_0(\lambda) d\lambda\right) = \int_{\lambda=-\infty}^{\infty} \frac{\partial^{k+l}[e^{i\lambda x - \lambda^2 t}]}{\partial t^k \partial x^l} \hat{u}_0(\lambda) d\lambda \quad (k, l \in \mathbb{N} \cup \{0\}).$$

Similar formulas hold also for the other integrals in the RHS of (4.1).

It follows from the above calculations that the function $u(x,t)$ is $C^\infty$ for $(x,t) \in Q$ and that it satisfies the equation $u_t = u_{xx}$, taking into consideration also the fact that the functions $e^{i\lambda x - \lambda^2 t}$, $\lambda \in \mathbb{C}$, satisfy this equation. Finally, differentiating (1.1) and (4.1), we obtain the formulas (1.2) and (1.3).

**Step 3** We will show that $\lim_{t \to 0^+} u(x,t) = u(x,t)\big|_{t=0} = u_0(x)$ ($x > 0$). (By writing «$u(x,t)\big|_{t=0}$» we mean the evaluation of the RHS of (1.2) at $t = 0$.) By Lebesgue's dominated convergence theorem,

$$\lim_{t \to 0^+} \frac{1}{2\pi} \int_{\lambda \in \Gamma} e^{i\lambda x - \lambda^2 t} \hat{u}_0(-\lambda) d\lambda = \frac{1}{2\pi} \int_{\lambda \in \Gamma} e^{i\lambda x} \hat{u}_0(-\lambda) d\lambda.$$

(This is justified by the fact that for $\lambda \in \Gamma$, $\left|e^{i\lambda x - \lambda^2 t}\right| = e^{-x \text{Im}\lambda} = e^{-x|\lambda|/\sqrt{2}}$.) Now by (3.10), the latter integral vanishes. Similarly, exploiting the presence of the factor $e^{i\lambda x}$, we obtain

$$\lim_{t \to 0^+} \frac{i}{\pi} \int_{\lambda \in \Gamma} e^{i\lambda x - \lambda^2 t} \lambda \tilde{g}_0(\lambda, t) d\lambda = 0.$$

But, by Fourier's inversion formula (2.1),

$$\lim_{t \to 0^+} \frac{1}{2\pi} \int_{\lambda=-\infty}^{\infty} e^{i\lambda x - \lambda^2 t} \hat{u}_0(\lambda) d\lambda = u_0(x) \text{ for } x > 0.$$

Also, it follows from Fourier's inversion formula (2.2) that

$$u(x,t)|_{t=0} = \left[\frac{1}{2\pi}\int_{\lambda=-\infty}^{\infty}e^{i\lambda x-\lambda^2 t}\hat{u}_0(\lambda)d\lambda\right]_{t=0} = \frac{1}{2\pi}\int_{\lambda=-\infty}^{\infty}e^{i\lambda x}\hat{u}_0(\lambda)d\lambda = u_0(x),$$

where the last integral is interpreted as the limit $\lim_{A\to\infty}\int_{-A}^{A}$.

Let us notice also that the RHS of (4.1), evaluated at $t=0$, is equal to $u_0(x)$ too. Indeed this follows from the fact that $\lambda\tilde{g}_0(\lambda,T) = \lambda\int_{\tau=0}^{T}e^{\lambda^2\tau}g_0(\tau)d\tau = O\left(\frac{1}{\lambda}\right)$ as $\lambda\to\infty$ with $\lambda\in\Omega^-$, which implies that

$$\left[\int_{\lambda\in\Gamma}e^{i\lambda x-\lambda^2 t}\lambda\tilde{g}_0(\lambda,T)d\lambda\right]_{t=0} = 0.$$

**Step 4** We claim that the RHS of (4.1), evaluated at $x=0$, is equal to $g_0(t)$, for $0<t<T$. Indeed, the aforementioned quantity, for $x=0$, becomes

$$\frac{1}{2\pi}\int_{\lambda=-\infty}^{\infty}e^{-\lambda^2 t}\hat{u}_0(\lambda)d\lambda - \frac{1}{2\pi}\int_{\lambda\in\Gamma}e^{-\lambda^2 t}\hat{u}_0(-\lambda)d\lambda - \frac{i}{\pi}\int_{\lambda\in\Gamma}e^{-\lambda^2 t}\lambda\tilde{g}_0(\lambda,T)d\lambda, \quad (5.4)$$

where the integrals over $\Gamma$ are interpreted as the corresponding limits $\lim_{A\to\infty}\int_{\lambda\in\Gamma_A}$. (In general, these two integrals do not converge absolutely since $\left|e^{-\lambda^2 t}\right|=1$ for $\lambda\in\Gamma$.)

The first two integrals in (5.4) are equal and therefore cancel each other, since by Cauchy's theorem in the set $\Omega^+\cap\{|\lambda|\leq A\}$,

$$\int_{-A}^{A}e^{-\lambda^2 t}\hat{u}_0(-\lambda)d\lambda - \int_{\lambda\in\Gamma_A}e^{-\lambda^2 t}\hat{u}_0(-\lambda)d\lambda = \int_{\Omega^+\cap K(0,A)}e^{-\lambda^2 t}\hat{u}_0(-\lambda)d\lambda,$$

and by Jordan's lemma 2.3.2,

$$\lim_{A\to\infty}\int_{\Omega^+\cap K(0,A)}e^{-\lambda^2 t}\hat{u}_0(-\lambda)d\lambda = 0.$$

(Notice that the value of the first integral in (5.4) does not change if we set $-\lambda$ in place of $\lambda$.)

It remains to show that the third integral in (5.4) is equal to $-g_0(t)$. Setting $\mu=i\lambda^2$, in this integral, we find

$$\frac{i}{\pi}\int_{\lambda\in\Gamma}e^{-\lambda^2 t}\lambda\tilde{g}_0(\lambda,T)d\lambda = \frac{1}{2\pi}\int_{\lambda\in\Gamma}e^{-\lambda^2 t}\int_{\tau=0}^{T}e^{\lambda^2\tau}g_0(\tau)d\tau(2i\lambda d\lambda)$$

$$= -\frac{1}{2\pi}\int_{\mu=-\infty}^{\infty}e^{i\mu t}\left(\int_{\tau=0}^{T}e^{-i\mu\tau}g_0(\tau)d\tau\right)d\mu = -g_0(t), \quad (5.5)$$

where the last equality follows from Fourier's inversion formula (2.2) applied to the function

$$f(\tau):=\begin{cases}g_0(\tau) & \text{for } 0\leq\tau\leq T \\ 0 & \text{for } \tau<0 \text{ or } \tau>T.\end{cases}$$

(The integral $\int_{\mu=-\infty}^{\infty}$ in (5.5) is interpreted as the limit $\lim_{A\to\infty}\int_{\mu=-A}^{A}$.) This completes the proof of the claim.

*Remark* The RHS of (1.1), evaluated at $x=0$, is equal to $g_0(t)/2$. Indeed,

$$\left[-\frac{i}{\pi}\int_{\lambda\in\Gamma}e^{i\lambda x-\lambda^2 t}\lambda\tilde{g}_0(\lambda,t)d\lambda\right]_{x=0} = -\frac{i}{\pi}\int_{\lambda\in\Gamma}e^{-\lambda^2 t}\lambda\tilde{g}_0(\lambda,t)d\lambda$$

$$= \frac{1}{2\pi} \int_{\mu=-\infty}^{\infty} e^{i\mu t}\left(\int_{\tau=0}^{t} e^{-i\mu\tau} g_0(\tau)d\tau\right) d\mu = \frac{1}{2} g_0(t),$$

where the last equality follows from Fourier's inversion formula (2.2) applied – this time – to the function

$$f(\tau) := \begin{cases} g_0(\tau) & \text{for } 0 \le \tau \le t \\ 0 & \text{for } \tau < 0 \text{ or } \tau > t. \end{cases}$$

**Step 5** We will show that, with the function $u(x,t)$ defined by (1.1) – equivalently by (4.1) – for $x > 0$ and $t > 0$, we have

$$\lim_{x \to 0^+} u(x,t) = g_0(t) \quad (t > 0). \tag{5.6}$$

Since for $x > 0$,

$$\int_{\lambda \in \Gamma} e^{i\lambda x - \lambda^2 t} \hat{u}_0(-\lambda) d\lambda = \int_{\lambda=-\infty}^{\infty} e^{i\lambda x - \lambda^2 t} \hat{u}_0(-\lambda) d\lambda,$$

Lebesgue's dominated convergence theorem gives

$$\lim_{x \to 0^+} \left[ \int_{\lambda=-\infty}^{\infty} e^{i\lambda x - \lambda^2 t} \hat{u}_0(\lambda) d\lambda - \int_{\lambda \in \Gamma} e^{i\lambda x - \lambda^2 t} \hat{u}_0(-\lambda) d\lambda \right] = \int_{\lambda=-\infty}^{\infty} e^{-\lambda^2 t} \hat{u}_0(\lambda) d\lambda - \int_{\lambda=-\infty}^{\infty} e^{-\lambda^2 t} \hat{u}_0(-\lambda) d\lambda = 0 \quad (t > 0).$$

Therefore, it suffices to show that

$$\lim_{x \to 0^+} \int_{\lambda \in \Gamma} e^{i\lambda x - \lambda^2 t} \lambda \tilde{g}_0(\lambda, T) d\lambda = \int_{\lambda \in \Gamma} e^{-\lambda^2 t} \lambda \tilde{g}_0(\lambda, T) d\lambda \quad (T > t), \tag{5.7}$$

since, by (5.5), the last integral is equal to $-\dfrac{\pi}{i} g_0(t)$.

Writing $\Gamma = \Gamma_1 + \Gamma_0$, with $\Gamma_1 = \Gamma \cap \{|\lambda| \ge 1\}$ and $\Gamma_0 = \Gamma \cap \{|\lambda| \le 1\}$, we see that (5.7) follows from

$$\lim_{x \to 0^+} \int_{\lambda \in \Gamma_1} e^{i\lambda x - \lambda^2 t} \lambda \tilde{g}_0(\lambda, T) d\lambda = \int_{\lambda \in \Gamma_1} e^{-\lambda^2 t} \lambda \tilde{g}_0(\lambda, T) d\lambda. \tag{5.8}$$

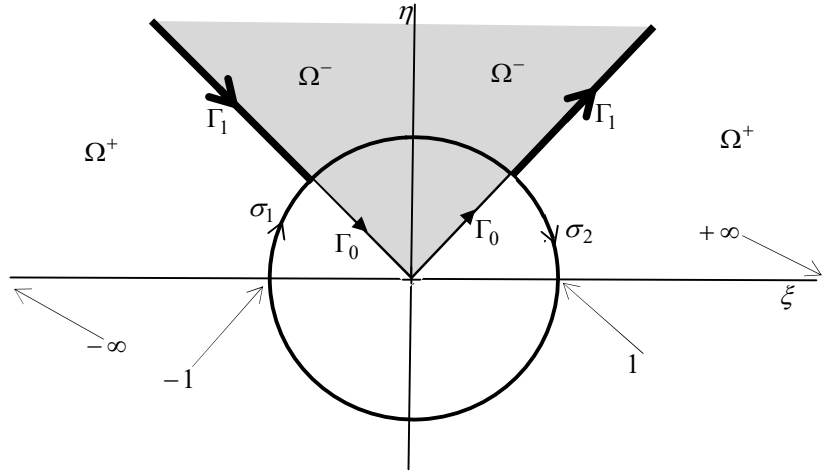

$\Gamma_1 = \Gamma \cap \{|\lambda| \ge 1\}$, $\Gamma_0 = \Gamma \cap \{|\lambda| \le 1\}$, $\sigma_1 = \{|\lambda| = 1\} \cap \Omega^+ \cap \{\text{Re } \lambda < 0\}$, $\sigma_2 = \{|\lambda| = 1\} \cap \Omega^+ \cap \{\text{Re } \lambda > 0\}$.

Integrating by parts we find that, for $\lambda \ne 0$ and $x \ge 0$,

$$e^{i\lambda x - \lambda^2 t} \lambda \tilde{g}_0(\lambda, T) = e^{i\lambda x - \lambda^2 t} \lambda \int_{\tau=0}^{T} e^{\lambda^2 \tau} g_0(\tau) d\tau$$

$$= \frac{1}{\lambda} e^{i\lambda x} e^{\lambda^2(T-t)} g_0(T) - \frac{1}{\lambda} e^{i\lambda x} e^{-\lambda^2 t} g_0(0) - e^{i\lambda x} \frac{1}{\lambda} e^{-\lambda^2 t} \int_{\tau=0}^{T} e^{\lambda^2 \tau} \frac{dg_0(\tau)}{d\tau} d\tau.$$

Integrating by parts once more and recalling that for $\lambda \in \Gamma_1$, $\left|e^{-\lambda^2 t}\right| = \left|e^{-\lambda^2 \tau}\right| = 1$ we find that

$$\frac{1}{\lambda} e^{-\lambda^2 t} \int_{\tau=0}^{T} e^{\lambda^2 \tau} \frac{dg_0(\tau)}{d\tau} d\tau = O(1/\lambda^3) \text{ for } \lambda \to \infty, \lambda \in \Gamma_1 \quad \text{(see also (2.10))}.$$

Therefore

$$\lim_{x \to 0^+} \int_{\Gamma_1} \left( e^{i\lambda x} \frac{1}{\lambda} e^{-\lambda^2 t} \int_{\tau=0}^{T} e^{\lambda^2 \tau} \frac{dg_0(\tau)}{d\tau} d\tau \right) d\lambda = \int_{\Gamma_1} \left( \frac{1}{\lambda} e^{-\lambda^2 t} \int_{\tau=0}^{T} e^{\lambda^2 \tau} \frac{dg_0(\tau)}{d\tau} d\tau \right) d\lambda.$$

Thus (5.8) will be proved if we show the following two equalities:

$$\lim_{x \to 0^+} \int_{\Gamma_1} \frac{1}{\lambda} e^{i\lambda x} e^{\lambda^2 (T-t)} d\lambda = \int_{\Gamma_1} \frac{1}{\lambda} e^{\lambda^2 (T-t)} d\lambda \text{ and } \lim_{x \to 0^+} \int_{\Gamma_1} \frac{1}{\lambda} e^{i\lambda x} e^{-\lambda^2 t} d\lambda = \int_{\Gamma_1} \frac{1}{\lambda} e^{-\lambda^2 t} d\lambda. \quad (5.9)$$

Proof of the 1$^{st}$ equality in (5.9). Since, by Cauchy's theorem in $\Omega^-$ and Jordan's lemma 2.3.1,

$$\int_{\Gamma_1} \frac{1}{\lambda} e^{i\lambda x} e^{\lambda^2 (T-t)} d\lambda - \int_{\{|\lambda|=1\} \cap \{\pi/4 \le \arg \lambda \le 3\pi/4\}} \frac{1}{\lambda} e^{i\lambda x} e^{\lambda^2 (T-t)} d\lambda = 0 \quad (x > 0),$$

it follows that

$$\lim_{x \to 0^+} \int_{\Gamma_1} \frac{1}{\lambda} e^{i\lambda x} e^{\lambda^2 (T-t)} d\lambda = \int_{\{|\lambda|=1\} \cap \{\pi/4 \le \arg \lambda \le 3\pi/4\}} \frac{1}{\lambda} e^{\lambda^2 (T-t)} d\lambda = \frac{1}{2} \int_{\{|\mu|=1\} \cap \{\operatorname{Im} \mu \le 0\}} \frac{1}{\mu} e^{-i\mu(T-t)} d\mu.$$

(The last equation follows by setting $\mu = i\lambda^2$.)

On the other hand,

$$\int_{\Gamma_1} \frac{1}{\lambda} e^{\lambda^2 (T-t)} d\lambda = -\frac{1}{2} \int_{-\infty}^{-1} \frac{1}{\mu} e^{-i\mu(T-t)} d\mu - \frac{1}{2} \int_{1}^{\infty} \frac{1}{\mu} e^{-i\mu(T-t)} d\mu.$$

(The above integrals do not converge absolutely and they have to be considered in the generalised sense.)

But, by Cauchy's theorem in the half-plane $\{z \in \mathbb{C} : \operatorname{Im} z \le 0\}$ and Jordan's lemma 2.3.1,

$$\int_{-\infty}^{-1} \frac{1}{\mu} e^{-i\mu(T-t)} d\mu + \int_{\{|\mu|=1\} \cap \{\operatorname{Im} \mu \le 0\}} \frac{1}{\mu} e^{-i\mu(T-t)} d\mu + \int_{1}^{\infty} \frac{1}{\mu} e^{-i\mu(T-t)} d\mu = 0,$$

and the 1$^{st}$ equality in (5.9) follows.

Proof of the 2$^{nd}$ equality in (5.9). We have

$$\int_{\Gamma_1} \frac{1}{\lambda} e^{i\lambda x} e^{-\lambda^2 t} d\lambda = \int_{\lambda=-\infty}^{-1} \frac{1}{\lambda} e^{i\lambda x} e^{-\lambda^2 t} d\lambda + \int_{\sigma_1} \frac{1}{\lambda} e^{i\lambda x} e^{-\lambda^2 t} d\lambda + \int_{\sigma_2} \frac{1}{\lambda} e^{i\lambda x} e^{-\lambda^2 t} d\lambda + \int_{\lambda=1}^{\infty} \frac{1}{\lambda} e^{i\lambda x} e^{-\lambda^2 t} d\lambda \quad (x > 0),$$

whence

$$\lim_{x \to 0^+} \int_{\Gamma_1} \frac{1}{\lambda} e^{i\lambda x} e^{-\lambda^2 t} d\lambda = \int_{-\infty}^{-1} \frac{1}{\lambda} e^{-\lambda^2 t} d\lambda + \int_{\sigma_1} \frac{1}{\lambda} e^{-\lambda^2 t} d\lambda + \int_{\sigma_2} \frac{1}{\lambda} e^{-\lambda^2 t} d\lambda + \int_{1}^{\infty} \frac{1}{\lambda} e^{-\lambda^2 t} d\lambda$$

$$= \int_{\sigma_1} \frac{1}{\lambda} e^{-\lambda^2 t} d\lambda + \int_{\sigma_2} \frac{1}{\lambda} e^{-\lambda^2 t} d\lambda.$$

Setting $\mu = i\lambda^2$ we find

$$\int_{\sigma_1} \frac{1}{\lambda} e^{-\lambda^2 t} d\lambda + \int_{\sigma_2} \frac{1}{\lambda} e^{-\lambda^2 t} d\lambda = \frac{1}{2} \int_{\tilde{\sigma}_1} \frac{1}{\mu} e^{i\mu t} d\mu + \frac{1}{2} \int_{\tilde{\sigma}_2} \frac{1}{\mu} e^{i\mu t} d\mu = -\frac{1}{2} \int_{\{|\mu|=1\} \cap \{\operatorname{Im} \mu \ge 0\}} \frac{1}{\mu} e^{i\mu t} d\mu,$$

where $-\tilde{\sigma}_1 = \{|\mu| = 1\} \cap \{0 \le \arg \mu \le \pi/2\}$ and $-\tilde{\sigma}_2 = \{|\mu| = 1\} \cap \{\pi/2 \le \arg \mu \le \pi\}$.

On the other hand,

$$\int_{\Gamma_1} \frac{1}{\lambda} e^{-\lambda^2 t} d\lambda = -\frac{1}{2} \int_{\mu=1}^{\infty} \frac{1}{\mu} e^{i\mu t} d\mu - \frac{1}{2} \int_{\mu=-\infty}^{-1} \frac{1}{\mu} e^{i\mu t} d\mu,$$

and the 2$^{nd}$ equality in (5.9) follows from the fact

$$\int_{\mu=-\infty}^{-1} \frac{1}{\mu} e^{i\mu t} d\mu - \int_{\{|\mu|=1\}\cap\{\text{Im}\,\mu\geq 0\}} \frac{1}{\mu} e^{i\mu t} d\mu + \int_{\mu=1}^{\infty} \frac{1}{\mu} e^{i\mu t} d\mu = 0.$$

The proof of (5.7) – and therefore that of (5.6) – is complete.

*Remark* The 1$^{st}$ equality in (5.9) does not hold when $t = T$. Indeed,

$$\lim_{x\to 0^+} \int_{\Gamma_1} \frac{1}{\lambda} e^{i\lambda x} d\lambda = \frac{i\pi}{2} \quad \text{while} \quad \int_{\Gamma_1} \frac{1}{\lambda} d\lambda = \lim_{A\to\infty} \int_{\Gamma_1 \cap \{|\lambda|\leq A\}} \frac{1}{\lambda} d\lambda = 0. \tag{5.10}$$

To prove the 1$^{st}$ equality in (5.10) we write

$$\int_{\Gamma_1} \frac{1}{\lambda} e^{i\lambda x} d\lambda = \int_{\lambda=-\infty}^{-1} \frac{1}{\lambda} e^{i\lambda x} d\lambda + \int_{\lambda=1}^{+\infty} \frac{1}{\lambda} e^{i\lambda x} d\lambda + \int_{\sigma_1} \frac{1}{\lambda} e^{i\lambda x} d\lambda + \int_{\sigma_2} \frac{1}{\lambda} e^{i\lambda x} d\lambda$$

$$= i\left[\int_{\lambda=-\infty}^{-1} \frac{\sin(\lambda x)}{\lambda} d\lambda + \int_{\lambda=1}^{+\infty} \frac{\sin(\lambda x)}{\lambda} d\lambda\right] + \int_{\sigma_1} \frac{1}{\lambda} e^{i\lambda x} d\lambda + \int_{\sigma_2} \frac{1}{\lambda} e^{i\lambda x} d\lambda$$

$$= i\left[\int_{\lambda=-\infty}^{-x} \frac{\sin\lambda}{\lambda} d\lambda + \int_{\lambda=x}^{+\infty} \frac{\sin\lambda}{\lambda} d\lambda\right] + \int_{\sigma_1} \frac{1}{\lambda} e^{i\lambda x} d\lambda + \int_{\sigma_2} \frac{1}{\lambda} e^{i\lambda x} d\lambda.$$

Letting $x \to 0^+$ we have

$$\lim_{x\to 0^+} \int_{\Gamma_1} \frac{1}{\lambda} e^{i\lambda x} d\lambda = i\left[\int_{\lambda=-\infty}^{0} \frac{\sin\lambda}{\lambda} d\lambda + \int_{\lambda=0}^{+\infty} \frac{\sin\lambda}{\lambda} d\lambda\right] + \int_{\sigma_1} \frac{1}{\lambda} d\lambda + \int_{\sigma_2} \frac{1}{\lambda} d\lambda = i\pi - \frac{i\pi}{2} = \frac{i\pi}{2}.$$

The proof of the 2$^{nd}$ equality in (5.10) is easier – it follows from the equation $(\log\lambda)' = 1/\lambda$, $\lambda \in \mathbb{C} - (-\infty, 0]$.

Thus, in general,

$$\lim_{x\to 0^+} \int_{\lambda\in\Gamma} e^{i\lambda x - \lambda^2 t} \lambda \widetilde{g}_0(\lambda, t) d\lambda \neq \int_{\lambda\in\Gamma} e^{-\lambda^2 t} \lambda \widetilde{g}_0(\lambda, t) d\lambda.$$

**Step 6** We will show that, for each fixed $t > 0$, the function $u(x,t)$ belongs to the space $C^\infty([0,\infty))$, with respect to the variable $x$. Firstly, noticing that

$$\frac{\partial^n}{\partial x^n}\left(\int_{\lambda=-\infty}^{\infty} e^{i\lambda x - \lambda^2 t} \hat{u}_0(\lambda) d\lambda\right) = i^n \int_{\lambda=-\infty}^{\infty} \lambda^n e^{i\lambda x - \lambda^2 t} \hat{u}_0(\lambda) d\lambda \quad (x > 0),$$

we have that

$$\lim_{x\to 0^+} \frac{\partial^n}{\partial x^n}\left(\int_{\lambda=-\infty}^{\infty} e^{i\lambda x - \lambda^2 t} \hat{u}_0(\lambda) d\lambda\right) = i^n \int_{\lambda=-\infty}^{\infty} \lambda^n e^{-\lambda^2 t} \hat{u}_0(\lambda) d\lambda, \tag{5.11}$$

and therefore the function

$$I_1(x,t) := \int_{\lambda=-\infty}^{\infty} e^{i\lambda x - \lambda^2 t} \hat{u}_0(\lambda) d\lambda$$

is $C^\infty$ up to the point $x = 0$, i.e., it belongs to the space $C^\infty([0,\infty))$ (with respect to the variable $x$). In the case of this integral – in treating it as a function of $x$ – important role is played by the presence of the factor $e^{-\lambda^2 t}$.

Similarly we treat the second integral

$$I_2(x,t) := \int_{\lambda\in\Gamma} e^{i\lambda x - \lambda^2 t} \hat{u}_0(-\lambda) d\lambda = \int_{\lambda=-\infty}^{\infty} e^{i\lambda x - \lambda^2 t} \hat{u}_0(-\lambda) d\lambda = \int_{\lambda=-\infty}^{\infty} e^{-i\lambda x - \lambda^2 t} \hat{u}_0(\lambda) d\lambda.$$

Now we consider the third integral

$$I_3(x,t) := \int_{\lambda \in \Gamma} e^{i\lambda x - \lambda^2 t} \lambda \tilde{g}_0(\lambda,t) d\lambda = \int_{\lambda \in \Gamma} e^{i\lambda x - \lambda^2 t} \lambda \tilde{g}_0(\lambda,\mathrm{T}) d\lambda, \text{ for } 0 < t < \mathrm{T} \text{ and } x > 0.$$

By the calculations that we made in step 5, it follows that the function $I_3(x,t)$ is continuous up to the point $x = 0$. Now we will prove that this function is $C^1$ up to the point $x = 0$, by proving that the limit

$$\lim_{x \to 0^+} \left[ \frac{\partial I_3(x,t)}{\partial x} \right] \text{ exists.} \tag{5.12}$$

By (2.7) we have
$$\frac{\partial I_3(x,t)}{\partial x} = \int_{\lambda \in \Gamma} i\lambda e^{i\lambda x - \lambda^2 t} \lambda \tilde{g}_0(\lambda,t) d\lambda$$

$$= ig_0(t) \int_{\lambda \in \Gamma} e^{i\lambda x} d\lambda - ig_0(0) \int_{\lambda \in \Gamma} e^{i\lambda x} e^{-\lambda^2 t} d\lambda - i \int_{\lambda \in \Gamma} e^{i\lambda x - \lambda^2 t} \int_{\tau=0}^{t} e^{\lambda^2 \tau} \frac{dg_0(\tau)}{d\tau} d\tau d\lambda. \tag{5.13}$$

(Here $x > 0$ and the presence in the above integrals of the factor $e^{i\lambda x}$ is crucial. As a matter of fact part of the difficulty is that this factor, as $x \to 0^+$, tends to 1.)

Since, for fixed $x > 0$, $d\left(\frac{1}{ix\lambda} e^{i\lambda x}\right)/d\lambda = e^{i\lambda x}$, we have

$$\int_{\lambda \in \Gamma} e^{i\lambda x} d\lambda = 0. \tag{5.14}$$

Also

$$\int_{\lambda \in \Gamma} e^{i\lambda x} e^{-\lambda^2 t} d\lambda = \int_{\lambda = -\infty}^{\infty} e^{i\lambda x} e^{-\lambda^2 t} d\lambda \tag{5.15}$$

and

$$e^{-\lambda^2 t} \int_{\tau=0}^{t} e^{\lambda^2 \tau} \frac{dg_0(\tau)}{d\tau} d\tau = O\left(\frac{1}{\lambda^2}\right) \text{ (as } \lambda \to \infty, \lambda \in \Gamma). \tag{5.16}$$

Substituting (5.14) and (5.15) in (5.13), we obtain

$$\frac{\partial I_3(x,t)}{\partial x} = -ig_0(0) \int_{\lambda = -\infty}^{\infty} e^{i\lambda x} e^{-\lambda^2 t} d\lambda - i \int_{\lambda \in \Gamma} e^{i\lambda x - \lambda^2 t} \int_{\tau=0}^{t} e^{\lambda^2 \tau} \frac{dg_0(\tau)}{d\tau} d\tau d\lambda, \text{ for } x > 0.$$

Thus, taking into consideration (5.16), we conclude that

$$\lim_{x \to 0^+} \left[ \frac{\partial I_3(x,t)}{\partial x} \right] = -ig_0(0) \int_{\lambda = -\infty}^{\infty} e^{-\lambda^2 t} d\lambda - i \int_{\lambda \in \Gamma} e^{-\lambda^2 t} \int_{\tau=0}^{t} e^{\lambda^2 \tau} \frac{dg_0(\tau)}{d\tau} d\tau d\lambda, \tag{5.17}$$

i.e., (5.12) holds.

Moreover, one more integration by parts (as in (2.7)) gives

$$\int_{\lambda \in \Gamma} e^{-\lambda^2 t} \int_{\tau=0}^{t} e^{\lambda^2 \tau} \frac{dg_0(\tau)}{d\tau} d\tau d\lambda = \int_{\lambda \in \Gamma_0} + \int_{\lambda \in \Gamma_1} = \int_{\lambda \in \Gamma_0} e^{-\lambda^2 t} \int_{\tau=0}^{t} e^{\lambda^2 \tau} \frac{dg_0(\tau)}{d\tau} d\tau d\lambda$$

$$+ ig_0(t) \int_{\lambda \in \Gamma_1} \frac{d\lambda}{\lambda^2} - ig_0(0) \int_{\lambda \in \Gamma_1} e^{-\lambda^2 t} \frac{d\lambda}{\lambda^2} - i \int_{\lambda \in \Gamma_1} e^{-\lambda^2 t} \frac{1}{\lambda^2} \int_{\tau=0}^{t} e^{\lambda^2 \tau} \frac{d^2 g_0(\tau)}{d\tau^2} d\tau d\lambda$$

which implies that the function

$$\frac{\partial}{\partial t} \left[ \lim_{x \to 0^+} \frac{\partial}{\partial x} \left( \int_{\lambda \in \Gamma_1} e^{i\lambda x - \lambda^2 t} \lambda \tilde{g}_0(\lambda,t) d\lambda \right) \right]$$

is continuous for $t > 0$. (For this conclusion we used also the fact that

$$\int_{\lambda \in \Gamma_1} e^{-\lambda^2 t} \frac{d\lambda}{\lambda^2} = \int_{\lambda \in \mathbb{R}, |\lambda| \geq 1} + \int_{\sigma_1 + \sigma_2} .)$$

Now it is easy to conclude that the function $g_1(t) = \lim_{x \to 0^+} [\partial u(x,t)/\partial x]$ is $C^1$ for $t > 0$.

To continue, in order to show that the limit $\lim_{x \to 0^+}\left[\partial^2 I_3(x,t)/\partial x^2\right]$ also exists, let us notice that

$$\frac{\partial^2}{\partial x^2}\left(\int_{\lambda \in \Gamma_1} e^{i\lambda x - \lambda^2 t}\lambda \widetilde{g}_0(\lambda,t)d\lambda\right)$$

$$= -g_0(t)\int_{\lambda \in \Gamma_1}\lambda e^{i\lambda x}d\lambda + g_0(0)\int_{\lambda \in \Gamma_1}\lambda e^{i\lambda x}e^{-\lambda^2 t}d\lambda + \int_{\lambda \in \Gamma_1} e^{i\lambda x - \lambda^2 t}\lambda \int_{\tau=0}^{t} e^{\lambda^2 \tau}\frac{dg_0(\tau)}{d\tau}d\tau d\lambda. \quad (5.18)$$

Observing that

$$\int_{\lambda \in \Gamma_1}\lambda e^{i\lambda x}d\lambda = -\int_{\lambda \in \Gamma_0}\lambda e^{i\lambda x}d\lambda, \quad \int_{\lambda \in \Gamma_1}\lambda e^{i\lambda x}e^{-\lambda^2 t}d\lambda = \int_{\lambda=-\infty}^{\infty}\lambda e^{i\lambda x}e^{-\lambda^2 t}d\lambda - \int_{\lambda \in \Gamma_0}\lambda e^{i\lambda x}e^{-\lambda^2 t}d\lambda \quad (5.19)$$

(their proofs are similar to the those of (5.14) and (5.15)) and taking into considarations the results of the computations made in step 5 and in particular the proof of (5.7) (in order to deal with the last term in (5.18)), we see that, indeed, the limit $\lim_{x \to 0^+}\left[\partial^2 I_3(x,t)/\partial x^2\right]$ exists.

Continuing in this way – inductively – we show that the limit

$$\lim_{x \to 0^+}\left[\frac{\partial^n I_3(x,t)}{\partial x^n}\right] \text{ exists, for every nonnegative integer } n,$$

and defines a $C^\infty$ function for $t > 0$. This completes the proof of the 3rd assertion of Theorem 1.1.

*Remark* More generally than (5.19), we have

$$\int_{\lambda \in \Gamma_1}\lambda^N e^{i\lambda x}d\lambda = -\int_{\lambda \in \Gamma_0}\lambda^N e^{i\lambda x}d\lambda \text{ and } \int_{\lambda \in \Gamma_1}\lambda^N e^{i\lambda x}e^{-\lambda^2 t}d\lambda = \int_{\lambda=-\infty}^{\infty}\lambda^N e^{i\lambda x}e^{-\lambda^2 t}d\lambda - \int_{\lambda \in \Gamma_0}\lambda^N e^{i\lambda x}e^{-\lambda^2 t}d\lambda \quad (N \in \mathbb{N}).$$

**Step 7** We will show that for every $x > 0$, the limit

$$\lim_{t \to 0^+}\frac{\partial^n u(x,t)}{\partial t^n} \text{ exists, for every nonnegative integer } n, \quad (5.20)$$

thus proving that $u(x,t) \in C^\infty([0,\infty))$ with respect to $t$.

For $n = 0$, (5.20) was proved in step 3. In the general case, using (4.1) for the definition of $u(x,t)$, it is easy to see that the existence of the limit in (5.20) – as far as the second and the third integral of (4.1) are concerned – is immediate due to the presence in these integrals of the factor $e^{i\lambda x}$, since the integration – in these two integrals – is carried out for $\lambda \in \Gamma$, where $\left|e^{i\lambda x}\right| = e^{-x\,\text{Im}\,\lambda} = e^{-x|\lambda|/\sqrt{2}}$.

Thus it remains to show that the limit

$$\lim_{t \to 0^+}\frac{\partial^n}{\partial t^n}\left[\int_{\lambda=-\infty}^{\infty}e^{i\lambda x - \lambda^2 t}\hat{u}_0(\lambda)d\lambda\right] \text{ exists for } n = 1,2,3,\ldots \quad (5.21)$$

For $t > 0$,

$$\frac{\partial}{\partial t}\left[\int_{\lambda=-\infty}^{\infty}e^{i\lambda x - \lambda^2 t}\hat{u}_0(\lambda)d\lambda\right] = -\int_{\lambda=-\infty}^{\infty}\lambda^2 e^{i\lambda x - \lambda^2 t}\hat{u}_0(\lambda)d\lambda.$$

Since $u_0(x) \in S([0,\infty))$, integration by parts gives

$$\lambda^2\hat{u}_0(\lambda) = \int_{y=0}^{\infty}\lambda^2 e^{-i\lambda y}u_0(y)dy = i\int_{y=0}^{\infty}\lambda u_0(y)\frac{\partial(e^{-i\lambda y})}{\partial y}dy$$

$$= -i\lambda u_0(0) + \int_{y=0}^{\infty}\frac{du_0(y)}{dy}\frac{\partial(e^{-i\lambda y})}{\partial y}dy = -i\lambda u_0(0) - u_0'(0) - \int_{y=0}^{\infty}u''_0(y)e^{-i\lambda y}dy$$

$$= -i\lambda u_0(0) - u_0'(0) - [u_0'']^{\wedge}(\lambda) \text{ (for } \lambda \in \mathbb{R}\text{)},$$

and therefore

$$\frac{\partial}{\partial t}\left[\int_{\lambda=-\infty}^{\infty} e^{i\lambda x-\lambda^2 t}\hat{u}_0(\lambda)d\lambda\right] = iu_0(0)\int_{\lambda=-\infty}^{\infty}\lambda e^{i\lambda x-\lambda^2 t}d\lambda + u'_0(0)\int_{\lambda=-\infty}^{\infty} e^{i\lambda x-\lambda^2 t}d\lambda + \int_{\lambda=-\infty}^{\infty} e^{i\lambda x-\lambda^2 t}[u_0'']^{\wedge}(\lambda)d\lambda$$

$$= iu_0(0)\int_{\lambda\in\Gamma}\lambda e^{i\lambda x-\lambda^2 t}d\lambda + u'_0(0)\int_{\lambda\in\Gamma} e^{i\lambda x-\lambda^2 t}d\lambda - \int_{\lambda=-\infty}^{\infty} e^{i\lambda x-\lambda^2 t}[u_0'']^{\wedge}(\lambda)d\lambda \quad (5.22)$$

where for the last equation we used Cauchy's theorem and Jordan's lemmas 2.3.5 and 2.3.6.
But by Fourier's inversion formula (2.1),

$$\lim_{t\to 0^+}\int_{\lambda=-\infty}^{\infty} e^{i\lambda x-\lambda^2 t}[u_0'']^{\wedge}(\lambda)d\lambda = 2\pi u_0''(x) \quad (x>0). \tag{5.23}$$

It follows from (5.22) and (5.23) that the limit

$$\lim_{t\to 0^+}\frac{\partial}{\partial t}\left[\int_{\lambda=-\infty}^{\infty} e^{i\lambda x-\lambda^2 t}\hat{u}_0(\lambda)d\lambda\right] \text{ exists.}$$

This proves (5.21) when $n=1$.
It is also clear from the above calculations that $u_1(x) = \lim_{t\to 0^+}[\partial u(x,t)/\partial t]$ is a $C^\infty$ function for $x>0$.
Proceeding in a similar way – inductively – we prove (5.20) and complete the proof of the 4$^{\text{th}}$ assertion of Theorem 1.1.

**Step 8** We will show that the function $u(x,t)$ is rapidly decreasing as $x\to +\infty$, uniformly for $t$ in compact subsets of $(0,\infty)$, i.e., given $m$, $n$ and $\beta > \alpha > 0$,

$$\lim_{x\to\infty}\left[x^m\frac{\partial^n u(x,t)}{\partial x^n}\right] = 0 \text{ uniformly for } \alpha\leq t\leq\beta. \tag{5.24}$$

Firstly, we deal with the integral $I_1(x,t) = \int_{\lambda=-\infty}^{\infty} e^{i\lambda x-\lambda^2 t}\hat{u}_0(\lambda)d\lambda$. We have

$$x^m\frac{\partial^n I_1(x,t)}{\partial x^n} = i^n x^m\int_{\lambda=-\infty}^{\infty}\lambda^n e^{i\lambda x}e^{-\lambda^2 t}\hat{u}_0(\lambda)d\lambda = \frac{i^{n-m-1}}{x}\int_{\lambda=-\infty}^{\infty}\frac{d^{m+1}(e^{i\lambda x})}{d\lambda^{m+1}}\left[\lambda^n e^{-\lambda^2 t}\hat{u}_0(\lambda)\right]d\lambda$$

$$= \frac{(-1)^{m+1}i^{n-m-1}}{x}\int_{\lambda=-\infty}^{\infty} e^{i\lambda x}\frac{d^{m+1}}{d\lambda^{m+1}}\left[\lambda^n e^{-\lambda^2 t}\hat{u}_0(\lambda)\right]d\lambda. \tag{5.25}$$

The last equation in (5.25) follows from integration by parts and the fact that the boundary terms vanish:

$$\left(\frac{d^{m+1-s}(e^{i\lambda x})}{d\lambda^{m+1-s}}\frac{d^s}{d\lambda^s}\left[\lambda^n e^{-\lambda^2 t}\hat{u}_0(\lambda)\right]\right)\bigg|_{\lambda=-\infty} = \left(\frac{d^{m+1-s}(e^{i\lambda x})}{d\lambda^{m+1-s}}e^{i\lambda x}\frac{d^s}{d\lambda^s}\left[\lambda^n e^{-\lambda^2 t}\hat{u}_0(\lambda)\right]\right)\bigg|_{\lambda=+\infty} = 0 \; (s\in\mathbb{N}\cup\{0\}).$$

Now let us observe that the derivative $d^{m+1}[\lambda^n e^{-\lambda^2 t}\hat{u}_0(-\lambda)]/d\lambda^{m+1}$ in the last integral in (5.25) is a finite linear combination of terms of the form

$$\lambda^{\ell_1} t^{\ell_2} e^{-\lambda^2 t}\frac{d^{\ell_3}\hat{u}_0(\lambda)}{d\lambda^{\ell_3}} \text{ with } \ell_1,\ell_2,\ell_3\in\mathbb{N}\cup\{0\},$$

and that, for $\alpha\leq t\leq\beta$,

$$\left|\lambda^{\ell_1}t^{\ell_2}e^{-\lambda^2 t}\frac{d^{\ell_3}\hat{u}_0(\lambda)}{d\lambda^{\ell_3}}\right| \leq |\lambda|^{\ell_1}\beta^{\ell_2}e^{-\lambda^2\alpha}\sup_{-\infty<\lambda<+\infty}\left|\frac{d^{\ell_3}\hat{u}_0(\lambda)}{d\lambda^{\ell_3}}\right| \text{ and } \sup_{-\infty<\lambda<+\infty}\left|\frac{d^{\ell_3}\hat{u}_0(\lambda)}{d\lambda^{\ell_3}}\right| < +\infty.$$

Therefore (5.25) implies that, for some $\ell\in\mathbb{N}$,

$$\left|x^m\frac{\partial^n I_1(x,t)}{\partial x^n}\right| \leq \frac{1}{x}\int_{\lambda=-\infty}^{\infty}|\lambda|^\ell e^{-\lambda^2\alpha}d\lambda, \text{ for every } \alpha\leq t\leq\beta \text{ and } x>0,$$

whence

$$\lim_{x \to \infty} \left[ x^m \frac{\partial^n I_1(x,t)}{\partial x^n} \right] = 0, \text{ uniformly for } \alpha \leq t \leq \beta. \tag{5.26}$$

Similarly

$$\lim_{x \to \infty} \left[ x^m \frac{\partial^n I_2(x,t)}{\partial x^n} \right] = 0, \text{ uniformly for } \alpha \leq t \leq \beta, \tag{5.27}$$

since

$$I_2(x,t) = \int_{\lambda \in \Gamma} e^{i\lambda x - \lambda^2 t} \hat{u}_0(-\lambda) d\lambda = \int_{\lambda = -\infty}^{\infty} e^{i\lambda x - \lambda^2 t} \hat{u}_0(-\lambda) d\lambda.$$

Next we consider the third integral $I_3(x,t) = \int_{\lambda \in \Gamma} e^{i\lambda x - \lambda^2 t} \lambda \tilde{g}_0(\lambda,t) d\lambda$ and we write

$$x^m \frac{\partial^n I_3(x,t)}{\partial x^n} = \frac{(-1)^{m+1} i^{n-m-1}}{x} \int_{\lambda \in \Gamma_0} \frac{d^{m+1}(e^{i\lambda x})}{d\lambda^{m+1}} e^{-\lambda^2 t} \lambda^{n+1} \tilde{g}_0(\lambda,t) d\lambda$$

$$+ i^n x^m \int_{\lambda \in \Gamma_1} e^{i\lambda x - \lambda^2 t} \lambda^{n+1} \tilde{g}_0(\lambda,t) d\lambda \tag{5.28}$$

where $\Gamma_0 = \Gamma \cap \{|\lambda| \leq 1\}$ and $\Gamma_1 = \Gamma \cap \{|\lambda| \geq 1\}$.

Since $\left|e^{-\lambda^2 t}\right| = 1$ and $\left|e^{i\lambda x}\right| = e^{-x|\lambda|/\sqrt{2}}$, for $\lambda \in \Gamma_1$, the integral over $\Gamma_1$ can be estimated as follows:

$$\left| i^n x^m \int_{\lambda \in \Gamma_1} e^{i\lambda x - \lambda^2 t} \lambda^{n+1} \tilde{g}_0(\lambda,t) d\lambda \right| \leq x^m \int_{\lambda \in \Gamma_1} e^{-x|\lambda|/\sqrt{2}} |\lambda|^{n+1} \left| e^{-\lambda^2 t} \int_{\tau=0}^{t} e^{\lambda^2 \tau} g_0(\tau) d\tau \right| d|\lambda|$$

$$\leq x^m \int_{\lambda \in \Gamma_1} e^{-x|\lambda|/\sqrt{2}} |\lambda|^{n+1} \int_{\tau=0}^{\beta} |g_0(\tau)| d|\lambda| = x^m e^{-x/(2\sqrt{2})} \left( \int_{\tau=0}^{\beta} |g_0(\tau)| \right) \int_{\lambda \in \Gamma_1} e^{-|\lambda|/(2\sqrt{2})} |\lambda|^{n+1} d|\lambda|,$$

for $x \geq 1$, $\alpha \leq t \leq \beta$. Therefore

$$\lim_{x \to \infty} \left[ x^m \int_{\lambda \in \Gamma_1} e^{i\lambda x - \lambda^2 t} \lambda^{n+1} \tilde{g}_0(\lambda,t) d\lambda \right] = 0, \text{ uniformly for } \alpha \leq t \leq \beta. \tag{5.29}$$

On the other hand, integrating by parts repeatedly, we see that the integral over $\Gamma_0$ can be estimated as follows:

$$\left| \int_{\lambda \in \Gamma_0} \frac{d^{m+1}(e^{i\lambda x})}{d\lambda^{m+1}} e^{-\lambda^2 t} \lambda^{n+1} \tilde{g}_0(\lambda,t) d\lambda \right| = \left| B(\lambda,x,t) \right|_{\lambda=e^{i\pi/4}} + \left| B(\lambda,x,t) \right|_{\lambda=e^{3i\pi/4}} +$$

$$\left| \int_{\lambda \in \Gamma_0} e^{i\lambda x} \frac{d^{m+1}}{d\lambda^{m+1}} \left[ e^{-\lambda^2 t} \lambda^{n+1} \tilde{g}_0(\lambda,t) \right] d\lambda \right| \tag{5.30}$$

where $B(\lambda,x,t)$ comes from the boundary terms of the integration by part processes and is a finite linear combination of terms of the form

$$e^{i\lambda x} x^{\ell_1} \lambda^{\ell_2} t^{\ell_3} e^{-\lambda^2 t} \int_{\tau=0}^{t} e^{\lambda^2 \tau} \tau^{\ell_4} g_0(\tau) d\tau \quad (\ell_1, \ell_2, \ell_3, \ell_4 \in \mathbb{N} \cup \{0\}).$$

Because of the presence of the factor $e^{i\lambda x}$, whose absolute value at $\lambda = e^{i\pi/4}$ (and at $\lambda = e^{i3\pi/4}$) is equal to $e^{-x/\sqrt{2}}$, it is easy to see that

$$\lim_{x \to \infty} \left[ \left| B(\lambda,x,t) \right|_{\lambda=e^{i\pi/4}} + \left| B(\lambda,x,t) \right|_{\lambda=e^{3i\pi/4}} \right] = 0, \text{ uniformly for } \alpha \leq t \leq \beta. \tag{5.31}$$

Also, since $|\lambda| \leq 1$ for $\lambda \in \Gamma_0$,

$$\frac{1}{x}\left|\int_{\lambda\in\Gamma_0} e^{i\lambda x}\frac{d^{m+1}}{d\lambda^{m+1}}\left[e^{-\lambda^2 t}\lambda^{n+1}\widetilde{g}_0(\lambda,t)\right]d\lambda\right|\leq\frac{1}{x},\ \text{uniformly for }\alpha\leq t\leq\beta. \quad (5.32)$$

Now

$$\lim_{x\to\infty}\left[x^m\frac{\partial^n I_3(x,t)}{\partial x^n}\right]=0,\ \text{uniformly for }\alpha\leq t\leq\beta, \quad (5.33)$$

follows from (5.32), (5.31), (5.30), (5.29) and (5.28).
Finally (5.24) follows from (5.26), (5.27) and (5.33).

Examining the results of the previous steps, we see that the proof of Theorem 1.1 is complete. □

## 6. More on the boundary behaviour of the solution

**Theorem 6.1** *With the assumptions as in Theorem 1.1, the function $u(x,t)$ defined by* (1.1) *satisfies the following:*

*1st The convergence $\lim_{t\to 0^+} u(x,t) = u_0(x)$ is uniform for $x$ in compact subsets of $(0,\infty)$.*

*2nd The limit condition $\lim_{\substack{(x,t)\to(x_0,0)\\(x,t)\in Q}} u(x,t) = u_0(x_0)$ for every $x_0 > 0$.*

*3rd The convergence $\lim_{\substack{(x,t)\to(x_0,0)\\(x,t)\in Q}} u(x,t) = u_0(x_0)$ is uniform for $x_0$ in compact subsets of $(0,\infty)$.*

**Proof.** Extending $u_0(y)$ also for $y<0$, by setting $u_0(y):=0$ for $y<0$, we may write the first integral in the RHS of (1.2) in the following way

$$\frac{1}{2\pi}\int_{\lambda=-\infty}^{\infty} e^{i\lambda x-\lambda^2 t}\hat{u}_0(\lambda)d\lambda = \frac{1}{2\pi}\int_{\lambda=-\infty}^{\infty} e^{i\lambda x-\lambda^2 t}\left(\int_{y=-\infty}^{\infty} e^{-i\lambda y}u_0(y)dy\right)d\lambda$$

$$= \frac{1}{2\pi}\int_{y=-\infty}^{\infty}\left(\int_{\lambda=-\infty}^{\infty} e^{i\lambda(x-y)}e^{-\lambda^2 t}d\lambda\right)u_0(y)dy = u_0 * \varphi_{\sqrt{t}}(x) \quad (6.1)$$

where

$$\varphi(x) = \frac{1}{2\sqrt{\pi}}e^{-x^2/4}\ \text{and}\ \varphi_{\sqrt{t}}(x) = \frac{1}{\sqrt{t}}\varphi(x/\sqrt{t}) = \frac{1}{2\sqrt{\pi t}}e^{-x^2/4t}.$$

(We point out that (6.1) holds for every $x\in\mathbb{R}$.) In the above calculation we used Fubini's theorem and the fact that

$$\frac{1}{2\pi}\int_{\lambda=-\infty}^{\infty} e^{i\lambda(x-y)}e^{-\lambda^2 t}d\lambda = \frac{1}{2\sqrt{\pi t}}e^{-(x-y)^2/4t},\ \text{for }x,y\in\mathbb{R}\text{ and }t>0.$$

Noticing that $\int_{x=-\infty}^{+\infty}\varphi(x)dx = 1$ and applying [13, Theorem 7.3] to (6.1), we obtain that

$$\lim_{t\to 0^+}\frac{1}{2\pi}\int_{\lambda=-\infty}^{\infty} e^{i\lambda x-\lambda^2 t}\hat{u}_0(\lambda)d\lambda = u_0(x),\ \text{uniformy for }x\text{ in compact subsets of }(0,\infty). \quad (6.2)$$

To deal with the second integral in the RHS of (1.1), let us recall that if $a>0$ then

$$\left|e^{i\lambda x}\right| = e^{-(\operatorname{Im}\lambda)x} = e^{-|\lambda|x/\sqrt{2}} \leq e^{-|\lambda|a/\sqrt{2}},\ \text{for }\lambda\in\Gamma\text{ and }x\geq a,\text{ and }\sup_{\lambda\in\Gamma}|\hat{u}_0(-\lambda)| < \infty.$$

Also $\left|e^{-\lambda^2 \tau}\right| = 1$ for $\lambda\in\Gamma$, whence

$$\left|e^{-\lambda^2 t} - 1\right| = \left|\int_{\tau=0}^{t} -\lambda^2 e^{-\lambda^2 \tau} d\tau\right| \leq t|\lambda|^2.$$

Therefore, for $x \geq a$,

$$\left|\int_{\lambda \in \Gamma} e^{i\lambda x - \lambda^2 t} \hat{u}_0(-\lambda) d\lambda - \int_{\lambda \in \Gamma} e^{i\lambda x} \hat{u}_0(-\lambda) d\lambda\right| = \left|\int_{\lambda \in \Gamma} e^{i\lambda x} \hat{u}_0(-\lambda)[e^{-\lambda^2 t} - 1] d\lambda\right|$$

$$\leq t \sup_{\lambda \in \Gamma} |\hat{u}_0(-\lambda)| \int_{\lambda \in \Gamma} e^{-|\lambda|a/\sqrt{2}} |\lambda|^2 d|\lambda|.$$

Recalling that $\int_{\lambda \in \Gamma} e^{i\lambda x} \hat{u}_0(-\lambda) d\lambda = 0$ (by (3.10)), we obtain that

$$\lim_{t \to 0^+} \frac{1}{2\pi} \int_{\lambda \in \Gamma} e^{i\lambda x - \lambda^2 t} \hat{u}_0(-\lambda) d\lambda = 0 \text{ uniformly for } x \geq a. \tag{6.3}$$

For the third integral in (1.2), we have

$$\left|\int_{\lambda \in \Gamma} e^{i\lambda x - \lambda^2 t} \lambda \tilde{g}_0(\lambda, t) d\lambda\right| = \left|\int_{\lambda \in \Gamma} e^{i\lambda x - \lambda^2 t} \lambda \left(\int_{\tau=0}^{t} e^{\lambda^2 \tau} g_0(\tau) d\tau\right) d\lambda\right| \leq \left(\int_{\tau=0}^{t} |g_0(\tau)| d\tau\right) \left(\int_{\lambda \in \Gamma} e^{-|\lambda|a/\sqrt{2}} |\lambda| d|\lambda|\right),$$

and therefore

$$\lim_{t \to 0^+} \int_{\lambda \in \Gamma} e^{i\lambda x - \lambda^2 t} \lambda \tilde{g}_0(\lambda, t) d\lambda = 0 \text{ uniformly for } x \geq a. \tag{6.4}$$

Now the 1st assertion of the theorem follows from (6.2), (6.3) and (6.4).

To prove the 2nd assertion let us fix $x_0 > 0$. Then, given any $\varepsilon > 0$, by the 1st assertion, there exists $\delta(\varepsilon) > 0$ so that $|u(x,t) - u_0(x)| < \varepsilon$ for $0 < t < \delta(\varepsilon)$ and $\frac{x_0}{2} \leq x \leq \frac{3x_0}{2}$. By making $\delta(\varepsilon)$ smaller if necessary, we may also achieve $|u_0(x) - u_0(x_0)| < \varepsilon$ for $|x - x_0| < \delta(\varepsilon)$. It follows that

$$|u(x,t) - u_0(x_0)| < 2\varepsilon \text{ for } |x - x_0| < \delta(\varepsilon) \text{ and } 0 < t < \delta(\varepsilon),$$

and this proves the 2nd assertion.
The proof of the 3rd assertion is similar. □

**Theorem 6.2** *With the assumptions as in Theorem 1.1, the function $u(x,t)$ defined by (1.2) satisfies the following:*

*1st The convergence $\lim_{x \to 0^+} u(x,t) = g_0(t)$ is uniform for $t$ in compact subsets of $(0, \infty)$.*

*2nd The limit condition $\lim_{\substack{(x,t) \to (0,t_0) \\ (x,t) \in \overline{Q}}} u(x,t) = g_0(t_0)$ for every $t_0 > 0$.*

*3rd The convergence $\lim_{\substack{(x,t) \to (0,t_0) \\ (x,t) \in \overline{Q}}} u(x,t) = g_0(t_0)$ is uniform for $t_0$ in compact subsets of $(0, \infty)$.*

**Proof** To prove the 1st assertion, we will follow the proof of (5.6) in step 5 of the proof of theorem 1.1, and, examining it more carefully we will see that (5.6) is actually uniform for $t$ in compact subsets of $(0, \infty)$. Firstly, for a fixed $a > 0$,

$$\sup_{t \geq a} \left|\int_{\lambda=-\infty}^{\infty} e^{i\lambda x - \lambda^2 t} \hat{u}_0(\lambda) d\lambda - \int_{\lambda=-\infty}^{\infty} e^{-\lambda^2 t} \hat{u}_0(\lambda) d\lambda\right| \leq \int_{\lambda=-\infty}^{\infty} \left|e^{i\lambda x} - 1\right| e^{-\lambda^2 a} |\hat{u}_0(\lambda)| d\lambda,$$

and therefore, by Lebesgue's dominated convergence theorem,

$$\lim_{x \to 0^+} \int_{\lambda=-\infty}^{\infty} e^{i\lambda x - \lambda^2 t} \hat{u}_0(\lambda) d\lambda = \int_{\lambda=-\infty}^{\infty} e^{-\lambda^2 t} \hat{u}_0(\lambda) d\lambda, \text{ uniformly for } t \geq a.$$

Working similarly with the second integral in (1.1), we conclude that

$$\lim_{x \to 0^+} \left[ \int_{\lambda=-\infty}^{\infty} e^{i\lambda x - \lambda^2 t} \hat{u}_0(\lambda) d\lambda - \int_{\lambda \in \Gamma} e^{i\lambda x - \lambda^2 t} \hat{u}_0(-\lambda) d\lambda \right] = 0, \text{ uniformly for } t \geq a.$$

(See also the corresponding calculation in step 5 of the proof of Theorem 1.1.)

Thus, fixing $T > a > 0$, it suffices to show that

$$\lim_{x \to 0^+} \int_{\lambda \in \Gamma} e^{i\lambda x - \lambda^2 t} \lambda \tilde{g}_0(\lambda, T) d\lambda = \int_{\lambda \in \Gamma} e^{-\lambda^2 t} \lambda \tilde{g}_0(\lambda, T) d\lambda \quad \text{uniformly for } t \text{ with } a \leq t < T. \tag{6.5}$$

Since $\int_{\lambda \in \Gamma} = \int_{\lambda \in \Gamma_0} + \int_{\lambda \in \Gamma_1}$, we have to show that

$$\lim_{x \to 0^+} \int_{\lambda \in \Gamma_0} e^{i\lambda x - \lambda^2 t} \lambda \tilde{g}_0(\lambda, T) d\lambda = \int_{\lambda \in \Gamma_0} e^{-\lambda^2 t} \lambda \tilde{g}_0(\lambda, T) d\lambda \quad \text{uniformly for } t \text{ with } a \leq t < T \tag{6.6}$$

and

$$\lim_{x \to 0^+} \int_{\lambda \in \Gamma_1} e^{i\lambda x - \lambda^2 t} \lambda \tilde{g}_0(\lambda, T) d\lambda = \int_{\lambda \in \Gamma_1} e^{-\lambda^2 t} \lambda \tilde{g}_0(\lambda, T) d\lambda \quad \text{uniformly for } t \text{ with } a \leq t < T. \tag{6.7}$$

Proof of (6.6). Since $\left| e^{-\lambda^2 t} \right| = 1$ for $\lambda \in \Gamma_0$, we have

$$\sup_{t < T} \left| \int_{\lambda \in \Gamma_0} e^{i\lambda x - \lambda^2 t} \lambda \tilde{g}_0(\lambda, T) d\lambda - \int_{\lambda \in \Gamma_0} e^{-\lambda^2 t} \lambda \tilde{g}_0(\lambda, T) d\lambda \right| \leq \int_{\lambda \in \Gamma_0} \left| e^{i\lambda x} - 1 \right| \left| \lambda \tilde{g}_0(\lambda, T) \right| d|\lambda|,$$

and (6.6) follows.

Proof of (6.7). Let us recall that

$$e^{i\lambda x - \lambda^2 t} \lambda \tilde{g}_0(\lambda, T) = \frac{1}{\lambda} e^{i\lambda x} e^{\lambda^2 (T-t)} g_0(T) - \frac{1}{\lambda} e^{i\lambda x} e^{-\lambda^2 t} g_0(0) - e^{i\lambda x} \frac{1}{\lambda} e^{-\lambda^2 t} \int_{\tau=0}^{T} e^{\lambda^2 \tau} \frac{dg_0(\tau)}{d\tau} d\tau \quad (\lambda \neq 0).$$

and

$$\frac{1}{\lambda} e^{-\lambda^2 t} \int_{\tau=0}^{T} e^{\lambda^2 \tau} \frac{dg_0(\tau)}{d\tau} d\tau$$

$$= \frac{1}{\lambda^3} e^{\lambda^2 (T-t)} \frac{dg_0(\tau)}{d\tau} \bigg|_{\tau=T} - \frac{1}{\lambda^3} e^{-\lambda^2 t} \frac{dg_0(\tau)}{d\tau} \bigg|_{\tau=0} - \frac{1}{\lambda^3} e^{-\lambda^2 t} \int_{\tau=0}^{T} e^{\lambda^2 \tau} \frac{d^2 g_0(\tau)}{d\tau^2} d\tau.$$

It follows from the above equation

$$\sup_{t<T} \left| \frac{1}{\lambda} e^{-\lambda^2 t} \int_{\tau=0}^{T} e^{\lambda^2 \tau} \frac{dg_0(\tau)}{d\tau} d\tau \right| \leq \frac{1}{|\lambda|^3} \left[ \left| \frac{dg_0(\tau)}{d\tau} \right|_{\tau=T} + \left| \frac{dg_0(\tau)}{d\tau} \right|_{\tau=0} + \int_{\tau=0}^{T} \left| \frac{d^2 g_0(\tau)}{d\tau^2} \right| d\tau \right] \text{ for } \lambda \in \Gamma_1.$$

The above estimate and the inequality

$$\left| \int_{\lambda \in \Gamma_1} e^{i\lambda x} \frac{1}{\lambda} e^{-\lambda^2 t} \int_{\tau=0}^{T} e^{\lambda^2 \tau} \frac{dg_0(\tau)}{d\tau} d\tau d\lambda - \int_{\lambda \in \Gamma_1} \frac{1}{\lambda} e^{-\lambda^2 t} \int_{\tau=0}^{T} e^{\lambda^2 \tau} \frac{dg_0(\tau)}{d\tau} d\tau d\lambda \right|$$

$$\leq \int_{\lambda \in \Gamma_1} \left| e^{i\lambda x} - 1 \right| \left| \frac{1}{\lambda} e^{-\lambda^2 t} \int_{\tau=0}^{T} e^{\lambda^2 \tau} \frac{dg_0(\tau)}{d\tau} d\tau \right| d|\lambda|$$

imply that

$$\lim_{x \to 0^+} \int_{\lambda \in \Gamma_1} e^{i\lambda x} \frac{1}{\lambda} e^{-\lambda^2 t} \int_{\tau=0}^{T} e^{\lambda^2 \tau} \frac{dg_0(\tau)}{d\tau} d\tau d\lambda$$

$$= \int_{\lambda \in \Gamma_1} \frac{1}{\lambda} e^{-\lambda^2 t} \int_{\tau=0}^{T} e^{\lambda^2 \tau} \frac{dg_0(\tau)}{d\tau} d\tau d\lambda \quad \text{uniformly for } t < T. \tag{6.8}$$

Also, for $x > 0$,

$$\int_{\Gamma_1} \frac{1}{\lambda} e^{i\lambda x} e^{\lambda^2 (T-t)} d\lambda = \int_{\sigma} \frac{1}{\lambda} e^{i\lambda x} e^{\lambda^2 (T-t)} d\lambda$$

(where we have set $\sigma = \{|\lambda|=1\} \cap \{\pi/4 \leq \arg\lambda \leq 3\pi/4\}$) and

$$\int_{\Gamma_1}\frac{1}{\lambda}e^{i\lambda x}e^{-\lambda^2 t}d\lambda = \int_{\lambda=-\infty}^{-1}\frac{1}{\lambda}e^{i\lambda x}e^{-\lambda^2 t}d\lambda + \int_{\sigma_1}\frac{1}{\lambda}e^{i\lambda x}e^{-\lambda^2 t}d\lambda + \int_{\sigma_2}\frac{1}{\lambda}e^{i\lambda x}e^{-\lambda^2 t}d\lambda + \int_{\lambda=1}^{\infty}\frac{1}{\lambda}e^{i\lambda x}e^{-\lambda^2 t}d\lambda. \quad (6.9)$$

But

$$\sup_{t<T}\left|\int_{\lambda\in\sigma}\frac{1}{\lambda}e^{i\lambda x}e^{\lambda^2(T-t)}d\lambda - \int_{\lambda\in\sigma}\frac{1}{\lambda}e^{\lambda^2(T-t)}d\lambda\right| = \sup_{t<T}\int_{\lambda\in\sigma}\left|e^{i\lambda x}-1\right|\left|e^{\lambda^2(T-t)}\right|d|\lambda|$$

$$\leq \sup_{t<T}\left[\sup_{\lambda\in\sigma}\left|e^{\lambda^2(T-t)}\right|\right]\int_{\lambda\in\sigma}\left|e^{i\lambda x}-1\right|d|\lambda| \quad (6.10)$$

which implies that the above quantity tends to zero, as $x \to 0^+$.
Similar calculations, based on (6.9), show also that

$$\lim_{x\to 0^+}\int_{\Gamma_1}\frac{1}{\lambda}e^{i\lambda x}e^{-\lambda^2 t}d\lambda = \int_{\Gamma_1}\frac{1}{\lambda}e^{-\lambda^2 t}d\lambda \quad \text{uniformly for } a \leq t < T. \quad (6.11)$$

(See also the corresponding calculation in step 5 of the proof of Theorem 1.1.)
Thus (6.7) follows from (6.8), (6.10) and (6.11), and this completes the proof of the 1st assertion. (We also used the result of step 4 of the proof of Theorem 1.1, where we showed that the quantity (5.4) is equal to $g_0(t)$.)
Finally the proof of the 2nd and 3rd assertion is similar to the proof of the corresponding assertions of Theorem 6.1. □

## 7. Boundary values of the derivatives of the solution

**Theorem 7.1** *With the assumptions as in Theorem 1.1, the function $u(x,t)$ defined by* (1.1) *satisfies the following:*

*1st The limit condition* $\lim_{t\to 0^+}\dfrac{\partial^n u(x,t)}{\partial x^n} = \dfrac{d^n u_0(x)}{dx^n}$, *with the convergence being uniform for $x$ in compact subsets of $(0,\infty)$.*

*2nd The limit condition* $\lim_{\substack{(x,t)\to(x_0,0)\\(x,t)\in Q}}\dfrac{\partial^n u(x,t)}{\partial x^n} = \dfrac{d^n u_0(x_0)}{dx^n}$, *with the convergence being uniform for $x_0$ in compact subsets of $(0,\infty)$.*

**Proof** Differentiating (1.1) with respect to $x$, we obtain

$$\frac{\partial u(x,t)}{\partial x} = \frac{1}{2\pi}\int_{\lambda=-\infty}^{\infty}i\lambda e^{i\lambda x-\lambda^2 t}\hat{u}_0(\lambda)d\lambda - \frac{1}{2\pi}\int_{\lambda\in\Gamma}i\lambda e^{i\lambda x-\lambda^2 t}\hat{u}_0(-\lambda)d\lambda - \frac{i}{\pi}\int_{\lambda\in\Gamma}i\lambda e^{i\lambda x-\lambda^2 t}\lambda\tilde{g}_0(\lambda,t)d\lambda. \quad (7.1)$$

For $\lambda \in \mathbb{R}$,

$$i\lambda\hat{u}_0(\lambda) = \int_{y=0}^{\infty}u_0(y)i\lambda e^{-i\lambda y}dy = -\int_{y=0}^{\infty}u_0(y)\frac{d(e^{-i\lambda y})}{dy}dy$$

$$= u_0(0) + \int_{y=0}^{\infty}\frac{du_0(y)}{dy}e^{-i\lambda y}dy = u_0(0) + \left(\frac{du_0(y)}{dy}\right)^{\wedge}(\lambda),$$

and therefore

$$\frac{1}{2\pi}\int_{\lambda=-\infty}^{\infty}i\lambda e^{i\lambda x-\lambda^2 t}\hat{u}_0(\lambda)d\lambda = \frac{1}{2\pi}u_0(0)\int_{\lambda=-\infty}^{\infty}e^{i\lambda x-\lambda^2 t}d\lambda + \frac{1}{2\pi}\int_{\lambda=-\infty}^{\infty}e^{i\lambda x-\lambda^2 t}\left(\frac{du_0(y)}{dy}\right)^{\wedge}(\lambda)d\lambda. \quad (7.2)$$

Andreas Chatziafratis

Since for $t > 0$, $\int_{\lambda=-\infty}^{\infty} e^{i\lambda x - \lambda^2 t} d\lambda = \sqrt{\dfrac{\pi}{t}} e^{-x^2/4t}$, we have

$$\lim_{t \to 0^+} \int_{\lambda \in \Gamma} e^{i\lambda x - \lambda^2 t} d\lambda = 0 \quad (x > 0). \tag{7.3}$$

Now, by the result of step 3 in the proof of Theorem 1.1, applied with the derivative $du_0(x)/dx$, we obtain

$$\lim_{t \to 0^+} \dfrac{1}{2\pi} \int_{\lambda=-\infty}^{\infty} e^{i\lambda x - \lambda^2 t} \left(\dfrac{du_0(y)}{dy}\right)^{\wedge}(\lambda) d\lambda = \dfrac{du_0(x)}{dx}. \tag{7.4}$$

Also for $\lambda \in \mathbb{C}$ with $\operatorname{Im} \lambda \geq 0$,

$$i\lambda \hat{u}_0(-\lambda) = \int_{y=0}^{\infty} u_0(y) i\lambda e^{i\lambda y} dy = \int_{y=0}^{\infty} u_0(y) \dfrac{d(e^{i\lambda y})}{dy} dy$$

$$= -u_0(0) - \int_{y=0}^{\infty} \dfrac{du_0(y)}{dy} e^{i\lambda y} dy = -u_0(0) - \left(\dfrac{du_0(y)}{dy}\right)^{\wedge}(-\lambda),$$

whence

$$\int_{\lambda \in \Gamma} i\lambda e^{i\lambda x - \lambda^2 t} \hat{u}_0(-\lambda) d\lambda = -u_0(0) \int_{\lambda \in \Gamma} e^{i\lambda x - \lambda^2 t} d\lambda - \int_{\lambda \in \Gamma} e^{i\lambda x - \lambda^2 t} \left(\dfrac{du_0(y)}{dy}\right)^{\wedge}(-\lambda) d\lambda.$$

Thus

$$\lim_{t \to 0^+} \int_{\lambda \in \Gamma} i\lambda e^{i\lambda x - \lambda^2 t} \hat{u}_0(-\lambda) d\lambda = -\lim_{t \to 0^+} \int_{\lambda \in \Gamma} e^{i\lambda x - \lambda^2 t} \left(\dfrac{du_0(y)}{dy}\right)^{\wedge}(-\lambda) d\lambda$$

$$= \int_{\lambda \in \Gamma} e^{i\lambda x} \left(\dfrac{du_0(y)}{dy}\right)^{\wedge}(-\lambda) d\lambda = 0 \tag{7.5}$$

where we used (7.3) and (3.10).

Finally, by Lebesque's dominated convergence theorem,

$$\lim_{t \to 0^+} \int_{\lambda \in \Gamma} i\lambda e^{i\lambda x - \lambda^2 t} \lambda \tilde{g}_0(\lambda, t) d\lambda = 0. \tag{7.6}$$

It follows from (7.6), (7.5) and (7.4), that

$$\lim_{t \to 0^+} \dfrac{\partial u(x,t)}{\partial x} = \dfrac{du_0(x)}{dx} \quad \text{for } x > 0. \tag{7.7}$$

Examining the above calculations and taking into consideration also the details of the proof of the 1st assertion of Theorem 6.1, we see that the convergence in (7.7) is actually uniform for $x$ in compact subsets of $(0, \infty)$. This proves the 1st assertion of the theorem for $n = 1$. Proceeding with induction on $n$, we can easily complete the proof of the 1st assertion. The 2nd assertion follows easily from the 1st one as in the case of Theorem 6.1. □

**Theorem 7.2** *With the assumptions as in Theorem 1.1, the function $u(x,t)$ defined by (1.1) satisfies the following:*

*1st The limit condition* $\lim\limits_{x \to 0^+} \dfrac{\partial^n u(x,t)}{\partial t^n} = \dfrac{d^n g_0(t)}{dt^n}$, *with the convergence being uniform for $t$ in compact subsets of $(0, \infty)$.*

*2nd The limit condition* $\lim\limits_{\substack{(x,t) \to (0,t_0) \\ (x,t) \in Q}} \dfrac{\partial^n u(x,t)}{\partial t^n} = \dfrac{d^n g_0(t_0)}{dt^n}$, *with the convergence being uniform for $t_0$ in compact subsets of $(0, \infty)$.*

**Proof** Differentiating (4.1) with respect to $t$, we have

$$\frac{\partial^n u(x,t)}{\partial t^n} = \frac{1}{2\pi}\int_{\lambda=-\infty}^{\infty}(-\lambda^2)^n e^{i\lambda x-\lambda^2 t}\hat{u}_0(\lambda)d\lambda - \frac{1}{2\pi}\int_{\lambda\in\Gamma}(-\lambda^2)^n e^{i\lambda x-\lambda^2 t}\hat{u}_0(-\lambda)d\lambda$$

$$-\frac{i}{\pi}\int_{\lambda\in\Gamma}(-\lambda^2)^n e^{i\lambda x-\lambda^2 t}\lambda\tilde{g}_0(\lambda,\mathrm{T})d\lambda, \text{ for } x>0, t<\mathrm{T}, \text{ and } n\in\mathbb{N}. \quad (7.8)$$

Firstly, for a fixed $a > 0$,

$$\sup_{t\geq a}\left|\int_{\lambda=-\infty}^{\infty}(-\lambda^2)^n e^{i\lambda x-\lambda^2 t}\hat{u}_0(\lambda)d\lambda - \int_{\lambda=-\infty}^{\infty}(-\lambda^2)^n e^{-\lambda^2 t}\hat{u}_0(\lambda)d\lambda\right| \leq \int_{\lambda=-\infty}^{\infty}\left|e^{i\lambda x}-1\right||\lambda|^{2n}e^{-\lambda^2 a}|\hat{u}_0(\lambda)|d\lambda,$$

and therefore, by Lebesgue's dominated convergence theorem,

$$\lim_{x\to 0^+}\int_{\lambda=-\infty}^{\infty}(-\lambda^2)^n e^{i\lambda x-\lambda^2 t}\hat{u}_0(\lambda)d\lambda = \int_{\lambda=-\infty}^{\infty}(-\lambda^2)^n e^{-\lambda^2 t}\hat{u}_0(\lambda)d\lambda, \text{ uniformly for } t\geq a.$$

On the other hand, since

$$\int_{\lambda\in\Gamma}(-\lambda^2)^n e^{i\lambda x-\lambda^2 t}\hat{u}_0(-\lambda)d\lambda = \int_{\lambda=-\infty}^{\infty}(-\lambda^2)^n e^{i\lambda x-\lambda^2 t}\hat{u}_0(-\lambda)d\lambda = \int_{\lambda=-\infty}^{\infty}(-\lambda^2)^n e^{-i\lambda x-\lambda^2 t}\hat{u}_0(\lambda)d\lambda,$$

working similarly with the second integral in (7.8), we conclude that

$$\lim_{x\to 0^+}\left[\int_{\lambda\in\Gamma}(-\lambda^2)^n e^{i\lambda x-\lambda^2 t}\hat{u}_0(-\lambda)d\lambda - \int_{\lambda=-\infty}^{\infty}(-\lambda^2)^n e^{-\lambda^2 t}\hat{u}_0(\lambda)d\lambda\right] = 0, \text{ uniformly for } t\geq a.$$

Therefore it suffices to show that

$$\lim_{x\to 0^+}\int_{\lambda\in\Gamma}(-\lambda^2)^n e^{i\lambda x-\lambda^2 t}\lambda\tilde{g}_0(\lambda,\mathrm{T})d\lambda = -\frac{\pi}{i}\frac{d^n g_0(t)}{dt^n}, \quad (7.9)$$

uniformly for $t$ in compact subsets of $(0,\infty)$.

Proof of (7.9) in the case $n=1$. Integrating by parts we find that, for every $\lambda\in\mathbb{C}$ and $x>0$,

$$-\lambda^2 e^{i\lambda x-\lambda^2 t}\lambda\tilde{g}_0(\lambda,\mathrm{T}) = -\lambda e^{i\lambda x}e^{\lambda^2(\mathrm{T}-t)}g_0(\mathrm{T}) + \lambda e^{i\lambda x}e^{-\lambda^2 t}g_0(0) + e^{i\lambda x}e^{-\lambda^2 t}\int_{\tau=0}^{\mathrm{T}}\lambda e^{\lambda^2\tau}\frac{dg_0(\tau)}{d\tau}d\tau,$$

and therefore

$$\int_{\lambda\in\Gamma}(-\lambda^2)e^{i\lambda x-\lambda^2 t}\lambda\tilde{g}_0(\lambda,\mathrm{T})d\lambda = -g_0(\mathrm{T})\int_{\lambda\in\Gamma}\lambda e^{i\lambda x}e^{\lambda^2(\mathrm{T}-t)}d\lambda + g_0(0)\int_{\lambda\in\Gamma}\lambda e^{i\lambda x-\lambda^2 t}d\lambda$$

$$+\int_{\lambda\in\Gamma}\left[\lambda e^{i\lambda x}e^{-\lambda^2 t}\int_{\tau=0}^{\mathrm{T}}e^{\lambda^2\tau}\frac{dg_0(\tau)}{d\tau}d\tau\right]d\lambda. \quad (7.10)$$

Now, by Cauchy's theorem and Jordan's lemma 2.3.7,

$$\int_{\lambda\in\Gamma}\lambda e^{i\lambda x}e^{\lambda^2(\mathrm{T}-t)}d\lambda = 0, \text{ for } x>0. \quad (7.11)$$

Similarly, by Jordan's lemmas 2.3.5 and 2.3.6,

$$\int_{\lambda\in\Gamma}\lambda e^{i\lambda x-\lambda^2 t}d\lambda = \int_{\lambda=-\infty}^{\infty}\lambda e^{i\lambda x-\lambda^2 t}d\lambda,$$

and therefore

$$\lim_{x\to 0^+}\int_{\lambda\in\Gamma}\lambda e^{i\lambda x-\lambda^2 t}d\lambda = 0, \text{ uniformly for } t\geq a. \quad (7.12)$$

But, from the proof of the 1st assertion of Theorem 6.2, applied with the function $dg_0(t)/dt$, we have

$$\lim_{x\to 0^+}\int_{\lambda\in\Gamma}\left[\lambda e^{i\lambda x}e^{-\lambda^2 t}\int_{\tau=0}^{\mathrm{T}}e^{\lambda^2\tau}\frac{dg_0(\tau)}{d\tau}d\tau\right]d\lambda = -\frac{\pi}{i}\frac{dg_0(t)}{dt}, \quad (7.13)$$

uniformly for $t$ in compact subsets of $(0,\infty)$.

Now (7.10), (7.11), (7.12), (7.13) imply (7.9) in the case $n=1$. The general case easily follows inductively. This completes the proof of the 1st assertion. The 2nd assertion follows as in the proof of the analogous parts of Theorem 6.1. □

**Theorem 7.3** *With the notation and the assumptions as in Theorem 1.1, we have the following:*
*1st The convergence*

$$u_n(x) = \lim_{t \to 0^+} \frac{\partial^n u(x,t)}{\partial t^n} \quad (n \in \mathbb{N} \cup \{0\})$$

*is uniform for $x$ in compact subsets of $(0, \infty)$.*

*Moreover $u_n(x) = \dfrac{d^{2n} u_0(x)}{dx^{2n}}$ and $\lim_{\substack{(x,t) \to (x_0, 0) \\ (x,t) \in Q}} \dfrac{\partial^n u(x,t)}{\partial t^n} = u_n(x_0)$, with the convergence being uniform for $x_0$ in compact subsets of $(0, \infty)$.*

*2nd The convergence*

$$g_n(t) = \lim_{x \to 0^+} \frac{\partial^n u(x,t)}{\partial x^n} \quad (n \in \mathbb{N} \cup \{0\})$$

*is uniform for $t$ in compact subsets of $(0, \infty)$.*

*Moreover $g_{2n}(t) = \dfrac{d^n g_0(t)}{dt^n}$ and $\lim_{\substack{(x,t) \to (0, t_0) \\ (x,t) \in Q}} \dfrac{\partial^n u(x,t)}{\partial x^n} = g_n(t_0)$, with the convergence being uniform for $t_0$ in compact subsets of $(0, \infty)$.*

**Proof** The 1st assertion follows from Theorem 7.1 since $\dfrac{\partial^n u(x,t)}{\partial t^n} = \dfrac{\partial^{2n} u(x,t)}{\partial x^{2n}}$. (For the last part we have to work as in Theorem 6.1.)

Also, examining the proof of the 3rd assertion of Theorem 1.1 given in step 6, we easily check that the convergence is actually uniform for $t$ in compact subsets of $(0, \infty)$, thus proving the 2nd assertion. □

Now we can prove Theorem 1.2.

***Proof of Theorem 1.2*** Since

$$\frac{\partial^{k+\ell} u(x,t)}{\partial x^k \partial t^\ell} = \frac{\partial^{k+2\ell} u(x,t)}{\partial x^{k+2\ell}}, \text{ for } x > 0 \text{ and } t > 0, \tag{7.14}$$

the conclusion follows from Theorems 7.1 and 7.3. We are also using the fact that the boundary values of the derivatives (7.14) are continuous on $\partial Q - \{(0,0)\} = \{(x,t) \in \mathbb{R}^2 : x = 0 \text{ or } t = 0\} - \{(0,0)\}$. More precisely, setting

$$\left.\frac{\partial^{k+\ell} u(x,t)}{\partial x^k \partial t^\ell}\right|_{x=0} := g_{k+2\ell}(t) \text{ (for } t > 0\text{)} \quad \text{and} \quad \left.\frac{\partial^{k+\ell} u(x,t)}{\partial x^k \partial t^\ell}\right|_{t=0} := u_{k+2\ell}(x) \text{ (for } x > 0\text{)},$$

we obtain the desired extensions. □

## 8. Proof of Theorem 1.3

Throughout this proof let us keep in mind that

$$u(x,t) = \frac{1}{2\pi} I_1(x,t) - \frac{1}{2\pi} I_2(x,t) - \frac{i}{\pi} I_3(x,t), \text{ for } x > 0 \text{ and } t > 0, \tag{8.1}$$

where

$$I_1(x,t) = \int_{\lambda=-\infty}^{\infty} e^{i\lambda x - \lambda^2 t} \hat{u}_0(\lambda) d\lambda = \int_{\lambda \in \mathbb{R}, |\lambda| \geq 1} e^{i\lambda x - \lambda^2 t} \hat{u}_0(\lambda) d\lambda + \int_{\lambda=-1}^{1} e^{i\lambda x - \lambda^2 t} \hat{u}_0(\lambda) d\lambda,$$

$$I_2(x,t) = \int_{\lambda=-\infty}^{\infty} e^{i\lambda x - \lambda^2 t} \hat{u}_0(-\lambda) d\lambda = \int_{\lambda=-\infty}^{\infty} e^{-i\lambda x - \lambda^2 t} \hat{u}_0(\lambda) d\lambda = \int_{\lambda \in \mathbb{R}, |\lambda| \geq 1} e^{-i\lambda x - \lambda^2 t} \hat{u}_0(\lambda) d\lambda + \int_{\lambda=-1}^{1} e^{-i\lambda x - \lambda^2 t} \hat{u}_0(\lambda) d\lambda,$$

and

$$I_3(x,t) = \int_{\lambda \in \Gamma} e^{i\lambda x - \lambda^2 t} \lambda \tilde{g}_0(\lambda, t) d\lambda = \int_{\Gamma_1} e^{i\lambda x - \lambda^2 t} \lambda \tilde{g}_0(\lambda, t) d\lambda + \int_{\lambda \in \Gamma_0} e^{i\lambda x - \lambda^2 t} \lambda \tilde{g}_0(\lambda, t) d\lambda.$$

**Step 1** Assuming $u_0(0) = g_0(0)$, we will prove that $\lim_{\overline{Q} \ni (x,t) \to (0,0)} u(x,t) = u_0(0)$, and for this it suffices to show that

$$\lim_{Q \ni (x,t) \to (0,0)} u(x,t) = \lim_{x \to 0^+} \left( \lim_{\substack{t \to 0^+ \\ x > 0}} u(x,t) \right) = \lim_{t \to 0^+} \left( \lim_{\substack{x \to 0^+ \\ t > 0}} u(x,t) \right). \quad (8.2)$$

Working for $x > 0$ and $t > 0$, and writing

$$\int_{\lambda \in \mathbb{R}, |\lambda| \geq 1} e^{i\lambda x - \lambda^2 t} \hat{u}_0(\lambda) d\lambda = \int_{\lambda \in \mathbb{R}, |\lambda| \geq 1} e^{i\lambda x - \lambda^2 t} \left[ \frac{u_0(0)}{i\lambda} + \frac{u_0'(0)}{(i\lambda)^2} + \frac{1}{(i\lambda)^2} \int_{y=0}^{\infty} e^{-i\lambda y} u_0''(y) dy \right] d\lambda$$

$$= u_0(0) \int_{|\lambda| \geq 1} e^{i\lambda x - \lambda^2 t} \frac{1}{i\lambda} d\lambda + \int_{|\lambda| \geq 1} e^{i\lambda x - \lambda^2 t} \frac{u_0'(0)}{(i\lambda)^2} d\lambda + \int_{|\lambda| \geq 1} e^{i\lambda x - \lambda^2 t} \frac{1}{(i\lambda)^2} \widehat{\left( \frac{d^2 u_0}{dy^2} \right)}(\lambda) d\lambda$$

and

$$\int_{\lambda \in \mathbb{R}, |\lambda| \geq 1} e^{-i\lambda x - \lambda^2 t} \hat{u}_0(\lambda) d\lambda = \int_{\lambda \in \mathbb{R}, |\lambda| \geq 1} e^{-i\lambda x - \lambda^2 t} \left[ \frac{u_0(0)}{i\lambda} + \frac{u_0'(0)}{(i\lambda)^2} + \frac{1}{(i\lambda)^2} \int_{y=0}^{\infty} e^{-i\lambda y} u_0''(y) dy \right] d\lambda$$

$$= u_0(0) \int_{|\lambda| \geq 1} e^{-i\lambda x - \lambda^2 t} \frac{1}{i\lambda} d\lambda + \int_{|\lambda| \geq 1} e^{-i\lambda x - \lambda^2 t} \frac{u_0'(0)}{(i\lambda)^2} d\lambda + \int_{|\lambda| \geq 1} e^{-i\lambda x - \lambda^2 t} \frac{1}{(i\lambda)^2} \widehat{\left( \frac{d^2 u_0}{dy^2} \right)}(\lambda) d\lambda,$$

we obtain

$$I_1(x,t) - I_2(x,t) = u_0(0) \int_{\lambda \in \mathbb{R}, |\lambda| \geq 1} 2i \sin(\lambda x) e^{-\lambda^2 t} \frac{1}{i\lambda} d\lambda$$

$$+ \int_{\lambda \in \mathbb{R}, |\lambda| \geq 1} 2i \sin(\lambda x) e^{-\lambda^2 t} \frac{1}{(i\lambda)^2} \widehat{\left( \frac{d^2 u_0}{dy^2} \right)}(\lambda) d\lambda + \int_{\lambda=-1}^{1} 2i \sin(\lambda x) e^{-\lambda^2 t} \hat{u}_0(\lambda) d\lambda,$$

where, for the last equation, we used also the fact that

$$\int_{\lambda \in \mathbb{R}, |\lambda| \geq 1} \sin(\lambda x) e^{-\lambda^2 t} \frac{1}{\lambda^2} d\lambda = 0.$$

It follows that

$$I_1(x,t) - I_2(x,t) \approx u_0(0) \int_{\lambda \in \mathbb{R}, |\lambda| \geq 1} 2i \sin(\lambda x) e^{-\lambda^2 t} \frac{1}{i\lambda} d\lambda. \quad (8.3)$$

*Definition of the symbol* $\approx$: For two functions $U(x,t)$ and $V(x,t)$, defined for $x > 0$ and $t > 0$, we will write $U(x,t) \approx V(x,t)$ if and only if

$$\lim_{Q \ni (x,t) \to (0,0)} [U(x,t) - V(x,t)] = \lim_{x \to 0^+} \left( \lim_{\substack{t \to 0^+ \\ x > 0}} [U(x,t) - V(x,t)] \right) = \lim_{t \to 0^+} \left( \lim_{\substack{x \to 0^+ \\ t > 0}} [U(x,t) - V(x,t)] \right),$$

provided that all the above limits exist.

On the other hand, for $x > 0$ and $t > 0$,

$$\int_{\Gamma_1} e^{i\lambda x - \lambda^2 t} \lambda \widetilde{g}_0(\lambda,t) d\lambda = \int_{\Gamma_1} e^{i\lambda x} \left[ \frac{g_0(t)}{\lambda} - \frac{g_0(0)}{\lambda} e^{-\lambda^2 t} - \frac{1}{\lambda} e^{-\lambda^2 t} \int_{\tau=0}^{t} e^{\lambda^2 \tau} g_0'(\tau) d\tau \right] d\lambda$$

and

$$\int_{\Gamma_1} e^{i\lambda x} \frac{g_0(t)}{\lambda} d\lambda = \int_{\substack{|\lambda| \geq 1 \\ \frac{3\pi}{4} \leq \arg \lambda \leq \frac{\pi}{4}}} e^{i\lambda x} \frac{g_0(t)}{\lambda} d\lambda \approx 0.$$

Also

$$\int_{\Gamma_1} e^{i\lambda x} \left[ \frac{1}{\lambda} e^{-\lambda^2 t} \int_{\tau=0}^{t} e^{\lambda^2 \tau} g_0'(\tau) d\tau \right] d\lambda \approx 0 \quad \text{and} \quad \int_{\lambda \in \Gamma_0} e^{i\lambda x - \lambda^2 t} \lambda \widetilde{g}_0(\lambda,t) d\lambda \approx 0.$$

Therefore

$$I_3(x,t) \approx -g_0(0) \int_{\Gamma_1} e^{i\lambda x - \lambda^2 t} \frac{1}{\lambda} d\lambda$$

$$= -g_0(0) \int_{\lambda \in \mathbb{R}, |\lambda| \geq 1} e^{i\lambda x - \lambda^2 t} \frac{1}{\lambda} d\lambda + g_0(0) \int_{\substack{|\lambda| = 1 \\ \arg \lambda \in [\frac{3\pi}{4}, \pi] \cup [0, \frac{\pi}{4}]}} e^{i\lambda x - \lambda^2 t} \frac{1}{\lambda} d\lambda,$$

whence

$$I_3(x,t) \approx -g_0(0) \int_{\lambda \in \mathbb{R}, |\lambda| \geq 1} i \sin(\lambda x) e^{-\lambda^2 t} \frac{1}{\lambda} d\lambda. \tag{8.4}$$

Since we assume $u_0(0) = g_0(0)$, (8.2) follows from (8.1), (8.3) and (8.4).

**Step 2** Assuming $u_0(0) = g_0(0)$, we will show that $\lim_{\overline{Q} \ni (x,t) \to (0,0)} \frac{\partial u(x,t)}{\partial x} = \frac{du_0(x)}{dx}\bigg|_{x=0}$, by proving that

$$\lim_{Q \ni (x,t) \to (0,0)} \frac{\partial u(x,t)}{\partial x} = \lim_{x \to 0^+} \left( \lim_{\substack{t \to 0^+ \\ x > 0}} \frac{\partial u(x,t)}{\partial x} \right) = \lim_{t \to 0^+} \left( \lim_{\substack{x \to 0^+ \\ t > 0}} \frac{\partial u(x,t)}{\partial x} \right). \tag{8.5}$$

For $x > 0$ and $t > 0$, we have

$$\frac{\partial u(x,t)}{\partial x} = \frac{1}{2\pi} \frac{\partial I_1(x,t)}{\partial x} - \frac{1}{2\pi} \frac{\partial I_2(x,t)}{\partial x} - \frac{i}{\pi} \frac{\partial I_3(x,t)}{\partial x}, \tag{8.6}$$

$$\frac{\partial I_1(x,t)}{\partial x} = \int_{\lambda \in \mathbb{R}, |\lambda| \geq 1} (i\lambda) e^{i\lambda x - \lambda^2 t} \hat{u}_0(\lambda) d\lambda + \int_{\lambda = -1}^{1} (i\lambda) e^{i\lambda x - \lambda^2 t} \hat{u}_0(\lambda) d\lambda,$$

$$\frac{\partial I_2(x,t)}{\partial x} = \int_{\lambda \in \mathbb{R}, |\lambda| \geq 1} (-i\lambda) e^{-i\lambda x - \lambda^2 t} \hat{u}_0(\lambda) d\lambda + \int_{\lambda = -1}^{1} (-i\lambda) e^{-i\lambda x - \lambda^2 t} \hat{u}_0(\lambda) d\lambda,$$

and

$$\frac{\partial I_3(x,t)}{\partial x} = \int_{\Gamma_1} (i\lambda) e^{i\lambda x - \lambda^2 t} \lambda \widetilde{g}_0(\lambda,t) d\lambda + \int_{\lambda \in \Gamma_0} (i\lambda) e^{i\lambda x - \lambda^2 t} \lambda \widetilde{g}_0(\lambda,t) d\lambda.$$

It follows that

$$\frac{\partial I_1(x,t)}{\partial x} \approx \int_{\lambda \in \mathbb{R}, |\lambda| \geq 1} (i\lambda) e^{i\lambda x - \lambda^2 t} \hat{u}_0(\lambda) d\lambda$$

$$\approx \int_{\lambda \in \mathbb{R}, |\lambda| \geq 1} e^{i\lambda x - \lambda^2 t} \left[ u_0(0) + \frac{u_0'(0)}{i\lambda} + \frac{u_0''(0)}{(i\lambda)^2} + \frac{1}{(i\lambda)^2} \int_{y=0}^{\infty} e^{-i\lambda y} u_0'''(y) dy \right] d\lambda$$

$$\approx u_0(0) \int_{\lambda\in\mathbb{R},|\lambda|\geq 1} e^{i\lambda x-\lambda^2 t} d\lambda + u_0'(0) \int_{\lambda\in\mathbb{R},|\lambda|\geq 1} e^{i\lambda x-\lambda^2 t} \frac{d\lambda}{i\lambda},$$

$$\frac{\partial I_2(x,t)}{\partial x} \approx -u_0(0) \int_{\lambda\in\mathbb{R},|\lambda|\geq 1} e^{-i\lambda x-\lambda^2 t} d\lambda - u_0'(0) \int_{\lambda\in\mathbb{R},|\lambda|\geq 1} e^{-i\lambda x-\lambda^2 t} \frac{d\lambda}{i\lambda},$$

and

$$\frac{\partial I_3(x,t)}{\partial x} = \int_{\Gamma_1} e^{i\lambda x} \left[ ig_0(t) - ig_0(0)e^{-\lambda^2 t} - ie^{-\lambda^2 t} \int_{\tau=0}^{t} e^{\lambda^2 \tau} g_0'(\tau) d\tau \right] d\lambda + \int_{\lambda\in\Gamma_0} (i\lambda) e^{i\lambda x-\lambda^2 t} \lambda \widetilde{g}_0(\lambda,t) d\lambda$$

$$\approx -ig_0(0) \int_{\lambda\in\mathbb{R},|\lambda|\geq 1} e^{-i\lambda x-\lambda^2 t} d\lambda \approx -ig_0(0) \int_{\lambda\in\mathbb{R},|\lambda|\geq 1} \cos(\lambda x) e^{-\lambda^2 t} d\lambda.$$

(For the last relation we used also the fact that, since $\int_{\Gamma} e^{i\lambda x} ig_0(t) d\lambda = 0$, we have

$$\int_{\Gamma_1} e^{i\lambda x} ig_0(t) d\lambda = -\int_{\Gamma_0} e^{i\lambda x} ig_0(t) d\lambda.)$$

Thus (8.6) gives

$$\frac{\partial u(x,t)}{\partial x} \approx \frac{1}{\pi} [u_0(0) - g_0(0)] \int_{\lambda\in\mathbb{R},|\lambda|\geq 1} \cos(\lambda x) e^{-\lambda^2 t} d\lambda,$$

and this, in view of the assumption $u_0(0) = g_0(0)$, implies (8.5).

***Step 3*** Assuming

$$u_0(0) = g_0(0) \text{ and } \left.\frac{d^2 u_0(x)}{dx^2}\right|_{x=0} = \left.\frac{dg_0(t)}{dt}\right|_{t=0}, \tag{8.7}$$

we will show that

$$\lim_{\overline{Q}\ni(x,t)\to(0,0)} \frac{\partial^2 u(x,t)}{\partial x^2} = \left.\frac{d^2 u_0(x)}{dx^2}\right|_{x=0}$$

by proving that

$$\lim_{Q\ni(x,t)\to(0,0)} \frac{\partial^2 u(x,t)}{\partial x^2} = \lim_{x\to 0^+} \left( \lim_{\substack{t\to 0^+ \\ x>0}} \frac{\partial u^2(x,t)}{\partial x^2} \right) = \lim_{t\to 0^+} \left( \lim_{\substack{x\to 0^+ \\ t>0}} \frac{\partial^2 u(x,t)}{\partial x^2} \right). \tag{8.8}$$

Working as in the previous steps we obtain

$$\frac{\partial^2 I_1(x,t)}{\partial x^2} \approx \int_{\lambda\in\mathbb{R},|\lambda|\geq 1} (i\lambda)^2 e^{i\lambda x-\lambda^2 t} \hat{u}_0(\lambda) d\lambda$$

$$\approx \int_{\lambda\in\mathbb{R},|\lambda|\geq 1} e^{i\lambda x-\lambda^2 t} \left[ (i\lambda) u_0(0) + u_0'(0) + \frac{u_0''(0)}{i\lambda} + \frac{1}{i\lambda} \int_{y=0}^{\infty} e^{-i\lambda y} u_0'''(y) dy \right] d\lambda$$

$$\approx u_0(0) \int_{\lambda\in\mathbb{R},|\lambda|\geq 1} (i\lambda) e^{i\lambda x-\lambda^2 t} d\lambda + u_0'(0) \int_{\lambda\in\mathbb{R},|\lambda|\geq 1} e^{i\lambda x-\lambda^2 t} d\lambda + u_0''(0) \int_{\lambda\in\mathbb{R},|\lambda|\geq 1} e^{i\lambda x-\lambda^2 t} \frac{d\lambda}{i\lambda},$$

$$\frac{\partial^2 I_2(x,t)}{\partial x^2} \approx \int_{\lambda\in\mathbb{R},|\lambda|\geq 1} (-i\lambda)^2 e^{-i\lambda x-\lambda^2 t} \hat{u}_0(\lambda) d\lambda$$

$$\approx \int_{\lambda\in\mathbb{R},|\lambda|\geq 1} e^{-i\lambda x-\lambda^2 t} \left[ (i\lambda) u_0(0) + u_0'(0) + \frac{u_0''(0)}{i\lambda} + \frac{1}{i\lambda} \int_{y=0}^{\infty} e^{-i\lambda y} u_0'''(y) dy \right] d\lambda$$

$$\approx u_0(0) \int_{\lambda\in\mathbb{R},|\lambda|\geq 1}(i\lambda)e^{-i\lambda x-\lambda^2 t}d\lambda + u_0'(0)\int_{\lambda\in\mathbb{R},|\lambda|\geq 1}e^{-i\lambda x-\lambda^2 t}d\lambda + u_0''(0)\int_{\lambda\in\mathbb{R},|\lambda|\geq 1}e^{-i\lambda x-\lambda^2 t}\frac{d\lambda}{i\lambda},$$

and

$$\frac{\partial^2 I_3(x,t)}{\partial x^2} \approx \int_{\Gamma_1}(i\lambda)^2 \lambda e^{i\lambda x}\left[\frac{g_0(t)}{\lambda^2}-\frac{g_0(0)}{\lambda^2}e^{-\lambda^2 t}-\frac{g_0'(t)}{\lambda^4}+\frac{g_0'(0)}{\lambda^4}e^{-\lambda^2 t}+\frac{1}{\lambda^4}e^{-\lambda^2 t}\int_{\tau=0}^{t}e^{\lambda^2\tau}g_0''(\tau)d\tau\right]d\lambda$$

$$\approx \int_{\Gamma_1}(i\lambda)^2\lambda e^{i\lambda x}\left[\frac{g_0(t)}{\lambda^2}-\frac{g_0(0)}{\lambda^2}e^{-\lambda^2 t}-\frac{g_0'(t)}{\lambda^4}+\frac{g_0'(0)}{\lambda^4}e^{-\lambda^2 t}\right]d\lambda$$

$$\approx -g_0(t)\int_{\Gamma_1}\lambda e^{i\lambda x}d\lambda + g_0(0)\int_{\Gamma_1}\lambda e^{i\lambda x-\lambda^2 t}d\lambda + g_0'(t)\int_{\Gamma_1}\frac{1}{\lambda}e^{i\lambda x}d\lambda - g_0'(0)\int_{\Gamma_1}\frac{1}{\lambda}e^{i\lambda x-\lambda^2 t}d\lambda$$

$$\approx g_0(0)\int_{\lambda\in\mathbb{R},|\lambda|\geq 1}\lambda e^{i\lambda x}e^{i\lambda x-\lambda^2 t}d\lambda - g_0'(0)\int_{\lambda\in\mathbb{R},|\lambda|\geq 1}\frac{1}{\lambda}e^{i\lambda x-\lambda^2 t}d\lambda,$$

where for the last relation we used also the following facts:

$$-g_0(t)\int_{\Gamma_1}\lambda e^{i\lambda x}d\lambda = g_0(t)\int_{\Gamma_0}\lambda e^{i\lambda x}d\lambda \approx 0 \text{ and } g_0'(t)\int_{\Gamma_1}\frac{1}{\lambda}e^{i\lambda x}d\lambda = g_0'(t)\int\frac{1}{\lambda}e^{i\lambda x}d\lambda \approx 0.$$

Therefore

$$\frac{\partial^2 I_1(x,t)}{\partial x^2}-\frac{\partial^2 I_2(x,t)}{\partial x^2}\approx -2u_0(0)\int_{\lambda\in\mathbb{R},|\lambda|\geq 1}\lambda\sin(\lambda x)e^{-\lambda^2 t}d\lambda + 2u_0''(0)\int_{\lambda\in\mathbb{R},|\lambda|\geq 1}\sin(\lambda x)e^{-\lambda^2 t}\frac{d\lambda}{\lambda}$$

and

$$\frac{\partial^2 I_3(x,t)}{\partial x^2}\approx ig_0(0)\int_{\lambda\in\mathbb{R},|\lambda|\geq 1}\lambda\sin(\lambda x)e^{-\lambda^2 t}d\lambda - ig_0'(0)\int_{\lambda\in\mathbb{R},|\lambda|\geq 1}\frac{1}{\lambda}\sin(\lambda x)e^{-\lambda^2 t}d\lambda,$$

whence (8.8) follows from (8.1) and (8.7).

**Step 4** In order to deal with the limit $\lim_{\bar{Q}\ni(x,t)\to(0,0)}\frac{\partial^3 u(x,t)}{\partial x^3}$, we compute

$$\frac{\partial^3 I_1(x,t)}{\partial x^3}\approx \int_{\lambda\in\mathbb{R},|\lambda|\geq 1}(i\lambda)^3 e^{i\lambda x-\lambda^2 t}\hat{u}_0(\lambda)d\lambda$$

$$\approx u_0(0)\int_{\lambda\in\mathbb{R},|\lambda|\geq 1}(i\lambda)^2 e^{i\lambda x-\lambda^2 t}d\lambda + u_0'(0)\int_{\lambda\in\mathbb{R},|\lambda|\geq 1}(i\lambda)e^{i\lambda x-\lambda^2 t}d\lambda$$

$$+u_0''(0)\int_{\lambda\in\mathbb{R},|\lambda|\geq 1}e^{i\lambda x-\lambda^2 t}d\lambda + u_0'''(0)\int_{\lambda\in\mathbb{R},|\lambda|\geq 1}e^{i\lambda x-\lambda^2 t}\frac{d\lambda}{i\lambda},$$

$$\frac{\partial^3 I_2(x,t)}{\partial x^3}\approx \int_{\lambda\in\mathbb{R},|\lambda|\geq 1}(-i\lambda)^3 e^{-i\lambda x-\lambda^2 t}\hat{u}_0(\lambda)d\lambda$$

$$\approx -u_0(0)\int_{\lambda\in\mathbb{R},|\lambda|\geq 1}(i\lambda)^2 e^{-i\lambda x-\lambda^2 t}d\lambda - u_0'(0)\int_{\lambda\in\mathbb{R},|\lambda|\geq 1}(i\lambda)e^{-i\lambda x-\lambda^2 t}d\lambda$$

$$-u_0''(0)\int_{\lambda\in\mathbb{R},|\lambda|\geq 1}e^{-i\lambda x-\lambda^2 t}d\lambda - u_0'''(0)\int_{\lambda\in\mathbb{R},|\lambda|\geq 1}e^{-i\lambda x-\lambda^2 t}\frac{d\lambda}{i\lambda}$$

and

$$\frac{\partial^3 I_3(x,t)}{\partial x^3}\approx \int_{\Gamma_1}(i\lambda)^3\lambda e^{i\lambda x}\left[\frac{g_0(t)}{\lambda^2}-\frac{g_0(0)}{\lambda^2}e^{-\lambda^2 t}-\frac{g_0'(t)}{\lambda^4}+\frac{g_0'(0)}{\lambda^4}e^{-\lambda^2 t}\right]d\lambda$$

$$\approx -ig_0(t)\int_{\Gamma_1}\lambda^2 e^{i\lambda x}d\lambda + ig_0(0)\int_{\Gamma_1}\lambda^2 e^{i\lambda x-\lambda^2 t}d\lambda + ig_0'(t)\int_{\Gamma_1}e^{i\lambda x}d\lambda - ig_0'(0)\int_{\Gamma_1}e^{i\lambda x-\lambda^2 t}d\lambda$$

$$\approx ig_0(0) \int_{\lambda\in\mathbb{R},|\lambda|\geq 1} \lambda^2 e^{i\lambda x-\lambda^2 t}d\lambda - ig_0'(0) \int_{\lambda\in\mathbb{R},|\lambda|\geq 1} e^{i\lambda x-\lambda^2 t}d\lambda.$$

Therefore

$$\frac{\partial^3 I_1(x,t)}{\partial x^3} - \frac{\partial^3 I_2(x,t)}{\partial x^3} \approx -2u_0(0) \int_{\lambda\in\mathbb{R},|\lambda|\geq 1} \lambda^2 \cos(\lambda x) e^{-\lambda^2 t}d\lambda + 2u_0''(0) \int_{\lambda\in\mathbb{R},|\lambda|\geq 1} \cos(\lambda x) e^{-\lambda^2 t}d\lambda$$

and

$$\frac{\partial^2 I_3(x,t)}{\partial x^2} \approx ig_0(0) \int_{\lambda\in\mathbb{R},|\lambda|\geq 1} \lambda^2 \cos(\lambda x) e^{-\lambda^2 t}d\lambda - ig_0'(0) \int_{\lambda\in\mathbb{R},|\lambda|\geq 1} \cos(\lambda x) e^{-\lambda^2 t}d\lambda,$$

and the equation

$$\lim_{\overline{Q}\ni(x,t)\to(0,0)} \frac{\partial^3 u(x,t)}{\partial x^3} = \left.\frac{d^3 u_0(x)}{dx^3}\right|_{x=0}$$

follows from (8.1) and (8.7). This completes the proof of Theorem 1.3. □

The above proof can easily be extended to prove the following generalization of Theorem 1.3.

**Theorem 8.1** *Let* $n \in \mathbb{N} \cup \{0\}$. *If* $u_0(x) \in \mathcal{S}([0,\infty))$ *and* $g_0(t) \in C^{\infty}([0,\infty))$ *such that*

$$\left.\frac{d^{2\ell}u_0(x)}{dx^{2\ell}}\right|_{x=0} = \left.\frac{d^{\ell}g_0(t)}{dt^{\ell}}\right|_{t=0} \text{ for } 0 \leq \ell \leq n,$$

*then the functions* $\partial^k u(x,t)/\partial x^k$ *(originally defined by (1.1) and (1.2) for* $(x,t)\in Q$ *and extended to* $\overline{Q}-\{(0,0)\}$ *by Theorem 1.2) satisfy the following:*

$$\lim_{\overline{Q}\ni(x,t)\to(0,0)} \frac{\partial^k u(x,t)}{\partial x^k} = \left.\frac{d^k u_0(x)}{dx^k}\right|_{x=0} \text{ for } 0 \leq k \leq 2n+1.$$

**Corolary 8.2** *With* $n$, $u_0$ *and* $g_0$ *as in the above theorem, the functions* $g_1(t)$, $g_3(t),...,$ $g_{2n+1}(t)$ *are continuous up to the point* $t=0$. *Also, if* $n \geq 1$,

$$\frac{dg_1(t)}{dt} = g_3(t), \frac{dg_3(t)}{dt} = g_5(t),..., \frac{dg_{2n-1}(t)}{dt} = g_{2n+1}(t),$$

*and* $g_{2n-1} \in C^1([0,\infty))$, $g_{2n-3} \in C^2([0,\infty)),..., g_1 \in C^n([0,\infty))$.

### Acknowledgements

The work contained in this preprint (January 2020) is part of my M.Sc. thesis (2018-19) at the National and Kapodistrian University of Athens, which has led to many novel results. I wish to express my gratitude to my Professors: N. Alikakos, G. Barbatis, T. Hatziafratis, I. Stratis, for insightful discussions and continuous support. The Onassis Foundation is also gratefully acknowledged for partial funding during that stage. The results presented herein have been published as:
A. Chatziafratis, D. Mantzavinos, Boundary behavior for the heat equation on the half-line, *Math. Methods Appl. Sci.* 45, 7364-93 (2022).